\documentclass[preprint,sort&compress,12pt]{elsarticle}
\usepackage{bbm}
\usepackage{amstext}

\usepackage{amsmath}
\usepackage{amssymb}
\usepackage{mathrsfs}
\usepackage{bm}
\usepackage{pstricks}
\usepackage{pst-coil}
\usepackage{pst-3d}
\usepackage{amsfonts}
\usepackage{amsthm}
\usepackage{CJK}
\usepackage{xcolor}
\allowdisplaybreaks
\newcommand{\mm}{\mathrm}

\newcommand{\mmd}{\mathrm{d}}

\newcommand{\be}{\begin{equation}}
\newcommand{\bea}{\begin{equation}\begin{aligned}}
\newcommand{\beas}{\begin{equation*}\begin{aligned}}
\newcommand{\eeas}{\end{aligned}\end{equation*}}
\newcommand{\eea}{\end{aligned}\end{equation}}
\newcommand{\ee}{\end{equation}}

\begin{document}
\begin{CJK*}{GBK}{song}
	\begin{frontmatter}
		\title{On the Inhibition of Rayleigh--Taylor Instability\\ by Capillarity in the Navier--Stokes--Korteweg Model}
		\author[FJ]{Fei Jiang}
		\ead{jiangfei0951@163.com}
		\author[FJ]{Yajie Zhang \corref{cor1}}
		\ead{zhangyajie315@qq.com}
		\author[SJ]{Zhipeng Zhang}\ead{zhangzhipeng@nju.edu.cn}
		\cortext[cor1]{Corresponding author. }
		\address[FJ]{School of Mathematics and Statistics, Fuzhou University, Fuzhou, 350108, China.}
		\address[SJ]{Department of Mathematics, Nanjing University, Nanjing 210093, China}
		\begin{abstract}			
Bresch--Desjardins--Gisclon--Sart had derived that the capillarity  slows  down the growth rate of Rayleigh--Taylor (RT) instability in an inhomogeneous
 incompressible fluid endowed with internal capillarity based on a \emph{linearized} incompressible Navier--Stokes--Korteweg (NSK) equations in 2008. Later Li--Zhang further obtained another result that the capillarity  inhibits RT instability also based on the \emph{linearized} equations in (SIAM J. Math. Anal. 3287--3315, 2023), if the capillarity coefficient is bigger than some threshold. In this paper, we further rigorously prove such phenomenon of capillarity inhibiting the RT instability  in the \emph{nonlinear}  incompressible NSK  equations in a horizontally periodic slab domain with Navier (slip) boundary conditions. The key idea in the proof is to capture the dissipative estimates of the tangential derivatives of density. Such dissipative estimates result in the decay-in-time of both the velocity and the perturbation density which is very useful to  overcome the difficulties arising from the nonlinear terms.
 \end{abstract}
	\end{frontmatter}
	\newtheorem{thm}{Theorem}[section]
	\newtheorem{lem}{Lemma}[section]
	\newtheorem{pro}{Proposition}[section]
	\newtheorem{concl}{Conclusion}[section]
	\newtheorem{cor}{Corollary}[section]
	\newproof{pf}{Proof}
	\newdefinition{rem}{Remark}[section]
	\newtheorem{definition}{Definition}[section]
	
\section{Introduction}\label{introud}
The equilibrium of a heavier fluid on the top of a lighter one, subject to gravity, is usually unstable. In fact, small disturbances acting on the equilibrium will grow and lead to the release of potential energy, as the heavier fluid moves down under gravity, and the lighter one is displaced upwards. This phenomenon was first studied by Rayleigh  \cite{RLISa} and then Taylor \cite{TGTP}, and is called therefore the Rayleigh--Taylor (RT) instability. In the last decades, the RT instability  {has been} extensively investigated from physical, numerical, and  mathematical aspects, see \cite{CSHHSCPO,WJH,desjardins2006nonlinear,hateau2005numerical,GBKJSAN} for examples and the references cited therein. It has been also widely analyzed {on how the} physical factors, such as elasticity \cite{JFJWGCOSdd}, rotation \cite{CSHHSCPO,BKASMMHRJA}, internal surface tension \cite{GYTI2,WYJTIKCT}, magnetic  {fields} \cite{JFJSWWWOA,JFJSJMFMOSERT,WYJ2019ARMA,JFJSARMA2019}, capillarity \cite{bresch2008instability} and so on, influence the dynamics of RT instability. In this paper, we are interested in the  phenomenon of capillarity inhibiting RT instability. Before stating our result and relevant progress in details, we need to mathematically formulate this  inhibition  phenomenon.

\subsection{Mathematical formulation for the capillary RT problem}\label{subsec:01}
\numberwithin{equation}{section}
A classical model to  describe the dynamics of an inhomogeneous
 incompressible fluid endowed with internal capillarity  (in the diffuse interface setting) in the presence of a uniform gravitational
field is the following  {general system of} incompressible Navier--Stokes--Korteweg (NSK) equations:
\begin{align}\label{0101}
\begin{cases}
\rho_t  +    \mathrm{div}(\rho v) = 0, \\
\rho v_t + \rho v\cdot \nabla v - \mathrm{div}(\mu (\rho )\mathbb{D}v)  + \nabla  {P}  = \mathrm{div}  {K} -\rho g\mathbf{e}^3,\\
\mathrm{div} v=0,
\end{cases}
\end{align}
where ${\rho (x,t)} \in \mathbb{R}^{+}$, $v(x,t) \in {\mathbb{R}^3}$ and $P(x,t)$ denote the density, velocity and kinetic pressure of the fluid {resp.} at the spacial position $x \in {\mathbb{R}^3}$ for time $t\in\mathbb{R}^+_0:=[0, {+\infty} )$. The differential operator $\mathbb{D}$ is defined by $\mathbb{D}v = \nabla v+\nabla v^{\top}$, where the subscript $\top$ denotes the transposition. $\mathbf{e}^3 $ represents the unit vector with the $3$-th component being $1$,  $g>0$ the gravitational constant  and $-\rho g \mathbf{e}^3$ the gravity. The shear viscosity function $\mu$ and the capillarity function $\kappa$ are known smooth functions $\mathbb{R}^+ \to \mathbb{R}$, and satisfy $\mu>0$  and $\kappa>0$. The general capillary tensor is written as
	\begin{align}
	K =  \left( {\rho \mm{div}(\kappa(\rho )\nabla \rho ) + \left( {\kappa(\rho ) - \rho\kappa'(\rho )} \right){{\left| {\nabla \rho } \right|}^2}} /2\right)\mathbb{I} -  {\kappa(\rho )\nabla \rho  \otimes \nabla \rho },\label{0101oooa}
	\end{align}
	where $\mathbb{I}$ denotes the identity matrix.    We mention that the well-posdeness problem for the incompressible NSK system has been investigated, see \cite{burtea2017lagrangian} and the references cited therein.
	
In classical hydrodynamics, the interface between two immiscible incompressible fluids is modeled as a free boundary which evolves in time. The equations describing the motion of each fluid are supplemented by boundary conditions at the free surface  involving the physical properties of the interface. For instance, in the free-boundary formulation, it is assumed that the interface has an internal surface tension. However, when the interfacial thickness is comparable to the length scale of the phenomena being examined, the free-boundary description breaks down. Diffuse-interface models provide an alternative description where the surface tension is expressed in its simplest form as $\mm{div}K$, i.e., the capillary tension which was introduced by  Korteweg in 1901 \cite{korteweg1901forme}. Later, its modern form was derived by Dunn and Serrin \cite{dunn1985thermomechanics}.  In addition, in the physical view, it can serve as a phase transition model to describe the motion of an incompressible fluid with capillarity effect.
	
To conveniently investigate the influence of capillarity on RT instability, we assume that \emph{$\mu$ and $\kappa$ are positive constants} as in \cite{bresch2008instability}, and thus get that
\begin{align}
 \mathrm{div}(\mu (\rho )\mathbb{D}v) =\mu\Delta v\mbox{ and } \mathrm{div}K =\kappa \rho\nabla \Delta \rho.
\end{align}
In addition, we consider the horizontally periodic   solutions of \eqref{0101}, and thus define a horizontally periodic domain via
\begin{align}\label{0101a}
\Omega:= 2\pi L_1\mathbb{T}\times 2\pi L_2\mathbb{T} \times(0,h),
\end{align}
where $\mathbb{T}:=\mathbb{R}/\mathbb{Z}$ and $L_i>0$ for $i=1$ and $2$. 	It is should be noted that if a function is defined on $\Omega$, then the function is horizontally periodic, i.e.,
	$$f(x_1,x_2,x_3)=f(2m\pi L_1+x_1, 2n\pi L_2+x_2,x_3)\mbox{ for any integers }m\mbox{ and }n.$$ \emph{We will see that the inhibiting effect of capillarity depends on
the periodic lengths $2\pi L_1$ and $2\pi L_2$} in Remark \ref{20224012302058}. The two-dimensional (2D) periodic domain
$2\pi L_1\mathbb{T}\times 2\pi L_2\mathbb{T}\times \{0,h\}$, denoted by $\partial\Omega$,  is customarily regarded as the  boundary of  the horizontally periodic domain $\Omega$.

We impose the following Navier (slip) boundary conditions for the velocity on $\partial\Omega$:
	\begin{gather}\label{n1c}
	v|_{\partial\Omega}\cdot{\mathbf{n}}=0,  \ ((\mathbb{D}v|_{\partial\Omega}){\mathbf{n}})_{\mm{tan}}=0,
	\end{gather}
	where  $\mathbf{n}$ denotes the outward unit normal vector to
	$\partial\Omega$, and the subscript ``$\mm{tan}$" means the tangential component of a vector  (for example,  $v_\mm{tan}=v-(v\cdot {\mathbf{n}}){\mathbf{n}}$) \cite{tapia2021stokes,ding2020rayleigh,ding2018stability,li2019global}.
	Since $\Omega$ is a slab domain, the Navier boundary condition is equivalent to the boundary condition
	\begin{align}
	(v_3,\partial_3 v_1,\partial_3 v_2)|_{\partial\Omega}=0.
	\label{20220202081737}
	\end{align}
We mention that the above boundary condition contributes to the mathematical verification of the inhibition phenomenon in this paper.

Let us further choose an  RT  equilibrium $(\bar{\rho},0)$ to \eqref{0101}, where the density profile $\bar{\rho}$ only depends on the third variable and satisfies
\begin{gather}
 \label{0102}
 \inf_{ x\in (0,h)}\{\bar{\rho}(x_3)\}>0,\\[1mm]
 \label{0102n}\bar{\rho}'|_{x_3=s}>0  \ \   \mbox{ for some} \ \ s\in (0,h).
\end{gather}
Here and in what follows $\bar{\rho}':=\mm{d}\bar{\rho}/\mm{d}x_3$. Then the pressure profile $\bar{P}$ under the equilibrium state is determined by the hydrostatic relation
\begin{equation}
	\nabla \bar{P}=\kappa\bar{\rho}\nabla \Delta \bar{\rho}-\bar{\rho}g  \mathbf{e}^3 \mbox{ in }  \Omega,  \label{equcomre}
	\end{equation}
which can be rewritten as an ordinary differential equation on $\bar{\rho}$:
$\bar{P}' =\kappa\bar{\rho}\bar{\rho}''' -\bar{\rho}g$.
We remark that the condition  \eqref{0102}
prevents us from treating vacuum, while the condition in \eqref{0102n} is called the RT condition,
which assures that there is at least a region in which the density is larger with increasing height $x_3$, thus leading to the classical RT instability, see \cite[Theorem 1.2]{JFJSO2014}.  \emph{However, we will see that such instability can be inhibited by capillarity} in Theorem \ref{thm2}.

Denoting the perturbations around the RT equilibrium by
 $$\varrho=\rho -\bar{\rho}, \ v= v- {0}, $$
 and recalling the relations of \eqref{equcomre} and
 \begin{align*} 	
 	  \nabla\rho\Delta \rho=  \nabla(\rho\Delta  \rho) - \rho\nabla \Delta  \rho,\end{align*}
 we obtain   the    perturbation system  from \eqref{0101}:
 \begin{equation}\label{1a}
 \begin{cases}
 \varrho _t+\bar{\rho}'v_3+v\cdot\nabla \varrho=0,\\
 (\bar{\rho}+{\varrho}) v_{t}+   (\bar{\rho}+{\varrho}) v\cdot  \nabla v+\nabla  \beta =\mu\Delta v-g{\varrho}\mathbf{e}^3-\kappa(\bar{\rho}''\nabla{\varrho}+\bar{\rho}'\Delta{\varrho}\mathbf{e}^3+\nabla{\varrho} \Delta{\varrho})
	, \\
\mm{div} v=0,
\end{cases}
\end{equation}
	where $\beta:= P+\kappa\bar{\rho}\Delta\bar{\rho}- \bar{P} -\kappa\rho\Delta \rho$.
	The corresponding initial and boundary conditions read as follows:
 \begin{align} \label{1ab}
	&(\varrho,v  )|_{t=0}=(\varrho^0,v^0 ) ,  \\
	\label{n1}
	& (v_3,\partial_3 v_1,\partial_3 v_2)|_{\partial\Omega}=0.
 \end{align}
Here and in what follows, we always use the right superscript $0$ to emphasize the initial data. We call the initial-boundary value problem \eqref{1a}--\eqref{n1} {\it the  capillary RT (abbr. CRT)   problem} for (the sake of) the simplicity. Obviously, to prove the phenomenon of the capillarity inhibiting the RT instability, it suffices to verify the stability in time of solutions to the above CRT problem
with non-trivial initial data.
	
In view of Li--Zhang's linear stability result for the CRT problem \emph{with a non-slip boundary condition} (in place of the Navier boundary condition) \cite{LFCZZP},  there exists  a threshold $\kappa_{\mm{C}}$ such that \emph{the linearized CRT system} (i.e. omitting the nonlinear perturbation terms in \eqref{1a}) is stable under the sharp stability condition
\begin{align}
\label{2020102241504} \kappa>\kappa_{\mm{C}},
\end{align} where we have defined that
\begin{align}
\label{2saf01504}
\kappa_{\mm{C}}:= {\sup_{ { w}\in H_{\sigma} }\frac{g\int\bar{\rho}' { w}_3^2 {\mm{d}x}}
{ \int |\bar{\rho}'\nabla  w_3|^2 {\mm{d}x}}}
\end{align}
		and $\bar{\rho}$ must satisfy the stabilizing condition
	\begin{align}
	\label{2022205071434}
	\inf_{ {x_3}\in (0,h)} \{|\bar{\rho}' {(x_3)}|\}>0.
	\end{align}
It should be noted that we have excluded the function $w$ satisfying $ \int |\bar{\rho}'\nabla  w_3  |^2 {\mm{d}x}=0$ in the above  definition of the supremum by default in \eqref{2saf01504}, and used the notations
\begin{align}
\int:=   \int_{(0,2\pi L_1)\times (0,2\pi L_2)\times (0,h)}\mbox{ and }H_{\sigma}:=\{w\in  {W^{1,2}(\Omega)}~|~\mm{div}w=0,\ w_3|_{\partial\Omega}=0\}.
\label{20224013022226}
\end{align}
However, if the following sharp instability condition is satisfied
  \begin{align}
\label{202010224150sdaf4} \kappa<\kappa_{\mm{C}},
\end{align}
there exists an unstable solution to the linearized CRT problem.

The linear stability result loosely presents that the capillarity can inhibit RT instability  in the fluid endowed with internal capillarity, if the capillarity coefficient is properly bigger; while
the linear instability result roughly shows that RT instability always occurs for too small capillarity coefficient. It is worth to point out that Bresch--Desjardins--Gisclon--Sart ever derived  that the capillarity  slows down the growth rate of linear RT instability in \cite{bresch2008instability}.

Compared with the linearized problem, more mathematical techniques shall be used to further analyze the instability/stability of the original nonlinear  CRT problem. For example, Nguyen has used the Guo--Strauss's bootstrap instability method in \cite{GYSWIC,GYSWICNonlinea} to establish the nonlinear RT instability result for the  CRT problem with \emph{non-slip boundary condition} under some sharp instability condition (similarly to \eqref{202010224150sdaf4}) \cite{nguyen2023influence}, also see the recent papers of Li--Zhang \cite{LFCZZP} and of Zhang \cite{zhang2022rayleigh}  for the nonlinear RT instability under the additional assumption $\kappa\ll 1$, and the recent papers of Zhang--Tian--Wang \cite{MR4645122}, resp. Zhang--Hua--Jiang--Lin \cite{MR4514792} for the nonlinear RT instability in NSK equations resp. Euler--Korteweg equations  with arbitrary value $\kappa>0$, where  the stabilizing condition \eqref{2022205071434} is violated.  Recently, motivated by the result of magnetic tension inhibiting the RT instability in the  2D  non-resistive magnetohydrodynamic fluid in \cite{JFJSZYYO}, Jiang--Li--Zhang \cite{jiang2023stability} mathematically proved that the capillarity  inhibits RT instability in the 2D NSK equations under the Lagrangian coordinates by making use of the dissipative estimates of $\eta_2$ and the near divergence-free condition of $\eta$ under small perturbation, where $\eta$ denotes the departure function of  fluid particles. However such method seems to be extremely difficulty to further verify the inhibition phenomenon for the 3D case, since we can not capture the dissipative estimates of $(\eta_1,\eta_2)$ from the dissipative estimates of $\eta_3$ and the near divergence-free condition of $\eta$. In this paper, we will develop a new mathematical  proof frame for the 3D case, and our stability result presents that the capillarity can inhibit RT instability in the fluid endowed with internal capillarity under the  sharp stability condition \eqref{2020102241504}  and the stabilizing condition	\eqref{2022205071434}, see Theorem \ref{thm2} for details.

 \subsection{Notations}\label{202402081729}

	Before stating our main result on  the CRT problem,  we shall introduce simplified notations which will be used throughout this paper.
	
	(1) Simplified basic notations:   $ \langle t\rangle:=1+t$,  $I_a:=(0,a)$ denotes a time interval, in particular, $I_\infty=\mathbb{R}^+:=(0,\infty)$.  $\overline{S}$ denotes the closure of a set $S\subset \mathbb{R}^n$ with $n\geqslant 1$, in particular, $\overline{I_T} =[0,T]$ and $\overline{I_\infty} = \mathbb{R}^+_0$.  $A:B:=a_{ij}b_{ij}$, where $A:=(a_{ij})_{n\times n}$, $B:=(b_{ij})_{n\times n}$ are $n\times n$ matrixes, and we have used the Einstein convention of summation over repeated indices.	$a\lesssim b$ means that $a\leqslant cb$ for some constant $c>0$, where $c >0$ may depend on  the domain $\Omega$, and the other known physical
functions/parameters, such as $\bar{\rho}$, $\mu$, $g$,  $\kappa$ in the CRT problem, and  vary from line to line.
	
$\partial_i:=\partial_{x_i}$, where $i=1$, $2$, $3$. Let $f:=(f_1,f_2,f_3)^{\top}$ be a vector function defined in a 3D domain, we define that $f_{\mm{h}}:=(f_1,f_2)^{\top}$ and
$\mm{curl}{f}:=(\partial_{2}f_3-\partial_{3}f_2,
\partial_{3}f_1-\partial_{1}f_3,\partial_{1}f_2-\partial_{2}f_1)^{\top}$. $\nabla_{\mm{h}}:=(\partial_{1},\partial_{2})^{\top}$. $\nabla^{\perp}_{\mm{h}}:=(-\partial_{2},\partial_{1})^{\top}$.  $\Delta_\mm{h}:=\partial^2_{1}+\partial^2_{2}$.   $\partial_{\mm{h}}^\alpha$ denotes $\partial_{1}^{\alpha_1}\partial_{2}^{\alpha_2}$
with $\alpha=(\alpha_1,\alpha_2)$, while $\partial_{\mm{h}}^i$ represents $\partial_{\mm{h}}^\alpha$ for any $|\alpha|:=\alpha_1+\alpha_2=i\in \mathbb{N}$. In particular $\partial_{\mm{h}}:=\partial^1_{\mm{h}}${.}	
	
$(v)_{\Omega}$ denotes the mean value of $v$ in a periodic cell $ (2\pi L_1\mathbb{T})\times( 2\pi L_2\mathbb{T})\times (0,h)$,  {i.e. $(v)_\Omega=\int  v\mathrm{d}{x}/{4\pi^2 L_1L_2h}$}.	 For the simplicity, we denote $\sqrt{\sum_{i=1}^n\|w_i\|_X^2}$ by $\|(w_1,\cdots,w_n)\|_X$, where $\|\cdot\|_X$ represents a norm or a semi-norm, and $w_i$ are scalar functions or functions for $1\leqslant i\leqslant n$.	
	
(2) Simplified Banach spaces, norms and semi-norms:
\begin{align}
&L^p:=L^p (\Omega)=W^{0,p}(\Omega),\
{H}^i:=W^{i,2}(\Omega ),   \nonumber \\[1mm]
&    {H}^j_{\mathrm{s}}:=\{w\in H^j_{\mm{s}}~|~ \partial_3w_1 =\partial_3w_2=w_3 =0\mbox{ on }\partial\Omega \},\ {_\sigma {H}^j_{\mathrm{s}}}:=H_\sigma\cap  {H}^j_{\mathrm{s}},\nonumber\\
&{ {H}}^4_{\bar{\rho}}:=\{\varrho\in H^4~|~\phi|_{\partial\Omega}=\partial^2_{3}\phi|_{\partial\Omega}=0,\ \partial_{3}\phi |_{\partial\Omega}=-\bar{\rho}'|_{\partial\Omega}\},\nonumber\\
& {^0_\sigma { {H}}^j_{\mathrm{s}}} : =\{w\in {_\sigma {H}^j_{\mathrm{s}}} ~|~({\rho}w_1)_\Omega=({\rho}w_2)_\Omega=0  \},\  \underline{H}^i=\{\phi\in H^i~|~(\phi)_\Omega=0\}, \nonumber  \\
&  \|\cdot \|_i :=\|\cdot \|_{H^i},\    \|\cdot\|_{{i},k}:=\sum_{\alpha_1+\alpha_2=i} \|\partial_{1}^{\alpha_1}\partial_{2}^{\alpha_2}\cdot\|_j , \nonumber
\end{align}
where $1\leqslant p\leqslant \infty$ is a real number, and $i\geqslant 0$, $j\geqslant 2$ are integers.

(3) Simplified spaces of functions with values in a Banach space:
\begin{align}
& L^p_TX:=L^p(I_T,X),\nonumber \\
&   {\mathfrak{P}} _{T}:=\{\varrho\in
C^0(\overline{I_T} ,   {H}^{4}  )~|~ \varrho_t\in C^0(\overline{I_T} ,H^3)\cap  L^2_TH^3\}, \nonumber\\
& {\mathcal{V}}_{ T}: =  \{v\in C^0(\overline{I_T},   {H}^3  )\cap L^2_T { {H}}^4  ~|~
v_t\in C^0(\overline{I_T} ,  {H}^1 )\cap  L^2_T {H}^2 \}. \nonumber
\end{align}
	
(4) Functionals of linearized potential energy: for $r\in H^1$,
\begin{align}&E_{\mm{L}}(r):=\kappa  \|\bar{\rho}'\nabla  r\|^2_0- g\int\bar{\rho}' r^2 \mm{d}x\label{edeE},\\
\label{eeE}
	&E(r):=\kappa\|\nabla r\|^2_{0}+\int\frac{\kappa\bar{\rho}'''-g}{\bar{\rho}'}r^2\mm{d}x .
	\end{align}
The function $E_{\mm{L}}(r)$ has been used by Li--Zhang for the proofs of the linear stablity/instably results \cite{LFCZZP}, and also by Jiang--Li--Zhang for the proofs of the nonlinear stablity/instably of the 2D NSK equations under the Lagrangian coordinates \cite{jiang2023stability}. However, we shall further use the new functional $E(r)$ for  the energy estimates and dissipative estimates of $\partial_{\mm{h}}\varrho$ (i.e. tangential derivatives of   density) for our 3D stability result in this paper.
		
(5) Energy and dissipation functionals (generalized):
\begin{align}
 \mathcal{E}:= \|\varrho\|^2_4+\|(\varrho_{t},v)\|^2_3 +\|v_t\|_1^2 \mbox{ and } \mathcal{D}:=\|v\|^2_4+\|\varrho\|^2_{1,3}+\|\varrho_t\|^2_{3}+\| v_t\|^2_2 .\label{di}			\end{align}
	
	(6) Both the functionals of tangential energy and tangential dissipation:
	\begin{align}
	 {\underline{\mathcal{E}}}:= \|\varrho\|_{{1},2}^2+\|\varrho_t\|_2^2+\|v_t\|_0^2+
\|(v_3,\partial_3 v_3)\|_1^2+   \|v\|_{{1},1}^2
\label{2022402200520957}
\end{align}
and
\begin{align}	\underline{\mathcal{D}} := \|\varrho\|_{1,1}^2+ \|\varrho\|_{{2},1}^2+ \|  (\varrho_t,v_t)  \|_1^2+ \|(v_2,\partial_3 v_3 )\|_{2}^2 +  \|v\|_{{1},2}^2   . \label{2022402011338}
	\end{align}
	
\subsection{A stability result}
Now  we state the stability result for the  CRT problem, which presents that the capillarity can inhibit the RT instability,	if the  capillary  coefficient  is properly large.
\begin{thm} \label{thm2}
Let $\mu$ and $\kappa$ be positive constants. If $\kappa$ and $\bar{\rho}\in {C^7}[0,h]$ satisfy \eqref{0102},  the sharp stability condition \eqref{2020102241504},
 the stabilizing condition \eqref{2022205071434}, and the additional boundary condition of the density profile
\begin{align}
\label{0102n1}	\bar{\rho}''|_{\partial\Omega}=0,
\end{align}
 there is a sufficiently small constant $\delta\in (0,1)$, such that for any $( \varrho^0,v^0)\in {H}^{4}_{\bar{\rho}} \times {^0_\sigma {{H}}^3_{\mathrm{s}}}$  satisfying a necessary compatibility condition	(i.e. $v_t(x,0)|_{\partial\Omega}=0$) and a smallness condition
\begin{align*}
\|(\nabla \varrho^0,v^0)\|_3 \leqslant\delta,
\end{align*}
the  CRT problem  of \eqref{1a}--\eqref{n1} admits a unique global classical  solution $(\varrho,v,\beta)\in {\mathfrak{P}} _{\infty}\times { \mathcal{V}_{\infty} }\times C^0(\mathbb{R}_0^+,\underline{H}^2)$. Moreover, the solution enjoys the stability estimate with algebraic decay-in-time
\begin{align}\label{1.200} \mathcal{E} (t)+ \langle t\rangle^2  {\underline{\mathcal{E}}} (t)+\int_0^t(\mathcal{D} (\tau)+ \langle \tau\rangle^2  \underline{\mathcal{D}} (\tau)   )\mm{d} \tau
\lesssim  \|(\nabla\varrho^0,v^0)\|_3^2\mbox{ for any }t>0.\end{align}
\end{thm}
\begin{rem}Thanks to the decay-in-time of dissipative estimates in \eqref{1.200}, we easily follow the argument of  (1.35) in \cite{JFJSZYYO} with slight modification to further derive the asymptotic behavior of the perturbation density:
\begin{align}
\sqrt{\langle t\rangle}\|\varrho(t)-\varrho^\infty\|_1\lesssim  \|(\nabla\varrho^0,v^0)\|_3 \mbox{ for some }\varrho^\infty\in H_{\bar{\rho}}^4 \mbox{ only  depending on }x_3.\label{1.200xx}
\end{align}
\end{rem}
\begin{rem}\label{20224012302058}
In view of Lemma \ref{261asdas567}, we easily find that
\begin{align}
\label{20223090318252}
0<\kappa_{\mm{C}}\leqslant  {g \|
\bar{\rho}'\|_{L^\infty}  \|
(\bar{\rho}')^{-1}\|_{L^\infty} ^2}(\pi^2h^{-2}+L_{\mm{max}}^{-2})^{-1} \mbox{ with  } L_{\mm{max}}:=\max\{L_1,L_2\}.
\end{align}
As a by-product, we observe that the smaller the  periodic cell is, the greater the stabilizing effect of capillarity is.  In particular, if the RT density profile is linear, then the threshold $\kappa_{\mm{C}}$ can be given by the formula $\kappa_{\mm{C}}=g /(\pi^2h^{-2}+L^{-2}_{\max})\bar{\rho}'$.
Thanks to such explicit expression of critical number, we easily guess that  the capillarity can also inhibit the RT instability in both the  infinite domains  $\mathbb{R}^2\times (0,h)$ (i.e. $L_{\mm{max}}=\infty$) and $2\pi L_1\mathbb{T}\times 2\pi L_2\mathbb{T} \times \mathbb{R}_+$  (i.e. $h=\infty$). However our proof for Theorem \ref{thm2} strongly depends on the finite periodic cell $(0,2\pi L_1)\times (0,2\pi L_2) \times(0,h)$, and thus can not be directly applied to  both the infinite domains. We shall develop a new proof frame for both the  infinite domains in an independent paper.	
\end{rem}	
\begin{rem}
It should be noted that  {in some situations $\kappa_{\mm{C}}=+\infty$. For example, we can choose $\bar{\rho}'\geqslant 0 $ such that }there exist four positive constants $\tilde{c}_1$, $\tilde{c}_2$, $s$, $\varepsilon$ and an interval   $(x_3^0-\varepsilon,x_3^0+\varepsilon)\subset (0,h)$ such that
$$\tilde{c}_1\leqslant  \frac{\bar\rho^\prime}{|x_3-x_3^0|^{2+s}}\leqslant \tilde{c}_2\mbox{
	for any }x_3\in  (x_3^0-\varepsilon,x_3^0+\varepsilon),$$
please refer to the proof of \cite[Proposition 2.1]{zhang2022rayleigh}. However, here  we exclude such case $\kappa_{\mm{C}}=+\infty$ by the stabilizing condition \eqref{2022205071434}, since we focus on the inhibiting effect of capillarity. In addition, if $\kappa>\kappa_{C}$  or $\kappa_{C}=\infty$,  we easily establish an RT  instability result   by following the arguments of Theorem 1.2 in \cite{JFJSZWC}  and Theorem 2.3 in \cite{nguyen2023influence}  under the absence of both the boundary conditions  \eqref{0102n1} and $\partial_{3}\varrho |_{\partial\Omega}=-\bar{\rho}'|_{\partial\Omega}$.
\end{rem}
\begin{rem}	
Both the boundary conditions  \eqref{0102n1} and $\partial_{3}\varrho |_{\partial\Omega}=-\bar{\rho}'|_{\partial\Omega}$ will be used to estimate for the highest-order spacial derivative of $(\varrho,v)$ in Lemma \ref{lem1}, which finally make sure us to further establish the dissipative estimate of $\|\varrho\|_{1,3}$ in  \eqref{r4}. Such dissipative estimate is extremely important to close the decay-in-time of  tangential estimates in \eqref{ed1A}.
\end{rem}

Now we  briefly sketch the proof of Theorem \ref{thm2}, the details of which will be
presented in Sections \ref{sec:global} and \ref{subsec:08a}.
Recalling the linear stability in \cite{LFCZZP}, we easily derive the basic  energy identity for the CRT problem defined on $I_T\times \Omega$ (see \eqref{2022401sdfa211727} for the derivation)
\begin{align}
 &\frac{1}{2}\frac{\mm{d}}{\mm{d}t}\left(E_{\mm{L}}\left(\int_0^tv_3(x,\tau)\mm{d}\tau\right)+\|\sqrt{\rho}  v\|^2_{0}  \right)  +\mu\int |\nabla v|^2\mm{d}x=\mathcal{I}_1(t), \end{align}
where we have defined that
\begin{align}	\mathcal{I}_1(t):=&  \int\left( \left(g\left(\int_0^tv\cdot \nabla \varrho \mm{d}s-\varrho^0\right)+\kappa\left( \bar{\rho}''' \left(\varrho^0-\int_0^tv\cdot \nabla \varrho \mm{d}s\right)\right.  \right.\right.
\nonumber\\
& \quad\left.\left.
\left. -\bar{\rho}'\Delta \left(\varrho^0-\int_0^tv\cdot \nabla \varrho \mm{d}s\right)\right)\right)\mathbf{e}^3
-\kappa\Delta{\varrho}\nabla{\varrho}
\right) \cdot v \mm{d}x.
\label{202402040453}
\end{align}
Thanks to the stabilizing condition \eqref{2020102241504} and the stabilizing condition	\eqref{2022205071434},   we easily deduce the following energy estimate from the above  energy identity:
\begin{align}	
&\|v\|_0^2+\left\| \int_0^tv_3(x,\tau)\mm{d}\tau \right\| ^2_1+c\int_0^t\| v(\tau)\|^2_1\mm{d} \tau  \lesssim \| \varrho^0\|_1^2+\|v^0\|_0^2 +\int_0^t \mathfrak{N}(\tau) \mm{d}\tau,\label{saf2022401211727}	\end{align}
where have defined that
\begin{align}	\mathfrak{N}(\tau) =
\|\varrho\|_{1,1}  \|\varrho\|_2\|v \|_1+ \left(\|\varrho\|_2^2+ \int_0^\tau
\left(\left\|\varrho\right\|_{1,1}\left\|v\right\|_2 +\left\|\varrho\right\|_3\left\|v_3\right\|_1 \right)(s)\mm{d}s   \right)\|v_3\|_1.\label{202402021099} \end{align}
In addition, by the mass equation \eqref{1a}$_1$ and the product estimate \eqref{fgestims}, we obtain that
\begin{align}
\|\varrho \|_0\lesssim \|\varrho^0\|_0+\left\|\int_0^tv_3(x,\tau)\mm{d}\tau \right\|_0 +\int_0^t ( \|\varrho\|_{1,1}\|v \|_1+\|\varrho\|_2 \|v_3\|_1) \mm{d}\tau,
\label{2022204022002023}
\end{align}
which, together with \eqref{saf2022401211727}, yields the basic (nonlinear) energy estimate:
 \begin{align}	
&\left\|\left(\varrho, v\right)\right\|_0^2+\left\|\int_0^tv_3(x,\tau)\mm{d}\tau\right\| ^2_1+c\int_0^t\| v(\tau)\|^2_1\mm{d} \tau
 \nonumber \\
 & \lesssim \| \varrho^0\|_1^2+\|v^0\|_0^2 +\left(\int_0^t ( \|\varrho\|_{1,1}\|v \|_1+\|\varrho\|_2 \|v_3\|_1) \mm{d}\tau\right)^2+	\int_0^t\mathfrak{N} (\tau)\mm{d}\tau .  \label{saf202sa2401211727}	\end{align}
Obviously, to control all the integrals from the nonlinear terms on the right hand  of the above inequality, it suffices to derive  the properly quick decay-in-time of dissipative estimates of $\partial_{\mm{h}}(\varrho,  v)$.

Before analyzing the decay-in-time, we shall first derive the dissipative estimate of $\partial_{\mm{h}}\varrho$. To this purpose,  we can apply curl operator to the momentum equation \eqref{1a}$_2$ and thus get the following vortex equation
\begin{equation}\label{nhp}
 \begin{aligned}
 &\kappa\bar{\rho}'\nabla_{\mm{h}}^{\perp}\Delta{\varrho}={\rho}\partial_{t} \omega_{\mm{h}} +(- \bar{\rho}'\partial_{t}v_1,  \bar{\rho}'\partial_{t}v_2 )^\top-\mu\Delta \omega_{\mm{h}}-(g-\kappa\bar{\rho}''')\nabla_{\mm{h}}^{\perp}{\varrho}+\rho v\cdot \nabla \omega_{\mm{h}}+ \mathbf{N}_{\mm{h}} ,
 \end{aligned}
 \end{equation}   see \eqref{202401120321011} for the definition of the nonlinear term  $\mathbf{N}$. Roughly speaking, the term $(g-\kappa\bar{\rho}''')\nabla_{\mm{h}}^{\perp}{\varrho}$ can be  controlled by stabilizing term $\kappa\bar{\rho}'\nabla_{\mm{h}}^{\perp}\Delta{\varrho}$ on the right hand of the above identity under the the stability condition \eqref{2020102241504} and   the stabilizing condition	\eqref{2022205071434} by using the new auxiliary functional $E $  in $\eqref{eeE}$. Thus the derivation of the dissipative estimates of $\partial_{\mm{h}}\varrho$ reduces to capture the dissipative estimates of $(v,v_t)$.
Fortunately we can build the dissipative estimates of $(v,v_t)$ under the Navier boundary condition.

By extremely refining the above analysis, we can use the energy method to arrive at the following total  energy inequality under small perturbation (referring to Proposition \ref{lem3}):
\begin{align}	
\sup_{0\leqslant  t\leqslant  T}{ \mathcal{E} (t)   }+\int_0^{T} \mathcal{D} (t) \mm{d}t
\lesssim& \|(\nabla\varrho^0,{v}^0)\|_3^2    + \sup_{0\leqslant  t\leqslant  T}\mathcal{E}(t)\int_0^{T}  \|v\|_{1,2} \mm{d}t \nonumber   \\
&  + \sup_{0\leqslant  t\leqslant  T}\sqrt{\mathcal{E}(t)}\left( \int_0^T ( \|\varrho\|_{1,1} +\|v_3\|_1 ) \mm{d}\tau\right)^2 .  \label{nhp1a}
\end{align}

Now we shall further establish the energy estimates of  $\partial_{\mm{h}}(\varrho,v)$ for the decay-in-time. To this end, we shall recall the other version of  basic energy identity  of the CRT problem (referring to \eqref{1e}):
\begin{align}	
\frac{\mm{d}}{\mm{d}t}(E(\varrho)+\|\sqrt{\rho}  v\|^2_{0})+2\mu\|\nabla v\|^2_{0} = \int \partial_3((\kappa\bar{\rho}'''-g)|\bar{\rho}'|^{-1}) \varrho^2 v_3\mm{d}x.\label{2022401211727}	\end{align}
Since the test function $w$ in the definition of $\kappa_{\mm{C}}$ should satisfy $\mm{div}w=0$, we can not apply the stabilizing condition \eqref{2020102241504} to $E(\varrho)$ for capture the energy estimate of $\varrho$. However, the stabilizing condition \eqref{2020102241504} can be applied to $E(\partial_{\mm{h}}\varrho)$. Based on this key observation and the tangential derivatives $\partial_{\mm{h}}(\varrho,v)$ enjoy an energy identity similarly to \eqref{2022401211727}, we can arrive at the decay-in-time of the tangential estimates under the small perturbation (referring Proposition \ref{lem2dx}):
\begin{align}
&\sup_{0\leqslant  t\leqslant  T}(\langle t\rangle^2 {	{\underline{\mathcal{E}}}(t)})+\int_0^{T}{\langle t\rangle^2 		\underline{\mathcal{D}} (t)} \mm{d}t\nonumber
\\
& \lesssim \|(\nabla\varrho^0,{v}^0)\|_3^2 + \int_{0}^{T} (\mathcal{D} (t)+
\sup_{0\leqslant  t\leqslant  T}(\langle t\rangle  \sqrt{\underline{\mathcal{E}}(t)})  \langle t\rangle \sqrt{\underline{\mathcal{D}}{\mathcal{D}}}))\mm{d}t,  \label{edds1A}
\end{align}
which, together with \eqref{nhp1a}, yields the desired stability estimate \eqref{1.200} under small perturbation.  In view of \eqref{1.200} and the unique local (-in-time) solvability of the CRT problem in Proposition \ref{202102182115}, we further obtain the unique global solvability with small perturbation immediately.

\section{\emph{A priori} estimates}\label{sec:global}
This section is devoted to  establishing both the total energy inequality   and the tangential energy inequality with decay-in-time  for the CRT problem \eqref{1a}--\eqref{n1}. To this purpose, let $(\varrho,v,q ) \in {\mathfrak{P}} _{T}\times { \mathcal{V}_{T} }\times C^0([0,T],\underline{H}^2)$   be a smooth solution to the  CRT  problem defined on  $I_T\times \Omega$.  We shall assume that the initial data $\varrho^0$ and $v^0$ belong  to ${ {\mathcal{H}}^4_{\bar{\rho}}}$ and ${^0_\sigma {H}^3_{\mathrm{s}}}$ resp.; moreover, $\varrho^0$  satisfies
\begin{equation}
0<{\inf\limits_{x_3\in(0,h))}\big\{\bar{\rho}(x_3) \big\}} \leqslant 2\inf\limits_{x\in\Omega}\big\{{\rho}^{0}(x)\big\} \leqslant 2\sup\limits_{x\in\Omega}\big\{{\rho}^{0}(x)\big\}  \leqslant 4\sup\limits_{x_3\in(0,h))}\big\{\bar{\rho}(x_3)\big\}  .\label{im1a}
\end{equation}
Here and in what follows $\rho^0:=\rho|_{t=0}$ and $\rho:= \bar{\rho}+\varrho$.
By the mass equation \eqref{1a}$_{1}$, it is well-known that
\begin{equation}
0<\inf\limits_{{x}\in\Omega}\big\{{\rho}^{0}({x})\big\}\leqslant {\rho}(t,x)\leqslant \sup\limits_{{x}\in\Omega}\big\{{\rho}^{0}({x})\big\}\mbox{ for any }(t,x)\in I_T\times\Omega.\label{im1}
\end{equation}
In addition, we shall keep in mind that $\mu$, $\kappa$ and $\bar{\rho}$ satisfy the assumptions of Theorem \ref{thm2}.

\subsection{Preliminaries}
In this section, we establish some preliminary results for  $(\varrho,v)$. To begin with, we derive more boundary conditions for $(\varrho, v_{\mm{h}},\omega)$.
\begin{lem}
	\label{2018051410720}
	The solution $(\varrho,{v})$ satisfies the following boundary conditions
	\begin{align}\label{varrho}
	&	(\varrho,\partial_3\rho, \partial_3^2\varrho)|_{\partial\Omega}=0,\\
	& \partial^3_{3}v_{\mm{h}}|_{\partial\Omega}=0,\label{2022401291254} \\
	&\label{omega} (\partial^{2i}_{3}\omega_{\mm{h}},\partial_{3}^{2i+1}\omega_3)|_{\partial\Omega}=0\mbox{ for }i=0,\ 1.
	\end{align}	
\end{lem}
\begin{pf}
	(1) We first derive \eqref{varrho}.  In view of the mass equation \eqref{1a}$_1$ and the boundary condition of $v_3$ in \eqref{n1}, it holds that
	\begin{equation}\label{n22a1}
	\varrho_t+v_{\mm{h}}\cdot \nabla_{\mm{h}}\varrho=0\  \text{on}\ \partial \Omega.
	\end{equation}
	Taking the inner product of the above identity
	and  $\varrho$ in $L^2(\partial\Omega)$, and then using the integration by parts and  the embedding inequality of $H^2\hookrightarrow C^0(\Omega)$ in \eqref{esmmdforinfty}, we derive that
	\begin{equation}\label{n22a2}
	\frac{\mm{d}}{\mm{d}t}	\int_{\partial\Omega}|\varrho|^2\mm{d}x_{\mm{h}}=-\frac{1}{2}\int_{\partial\Omega} v_{\mm{h}}\nabla_{\mathrm{h}} |\varrho|^2 \mm{d}x_{\mm{h}}  \lesssim \|\mm{div}_{\mathrm{h}}v_{\mm{h}}\|_{2} \int_{\partial\Omega} |\varrho|^2\mm{d}x_{\mm{h}}  .
	\end{equation}
	Noting that  $\int_0^T \|\mm{div}_{\mathrm{h}}v_{\mm{h}}\|_{2}\mm{d}\tau<\infty$ and $ \varrho^0|_{\partial\Omega} =0$, thus applying Growall's inequality to the above inequality yields
	\begin{equation}\label{n22asf}
	\|\varrho\|^2_{L^2(\partial\Omega)}=0,
	\end{equation}
	which implies
	\begin{equation}\label{n22a}
	\varrho|_{\partial\Omega}=0.
	\end{equation}
	
	Applying $\partial_3$ to \eqref{1a}$_1$, and then using  the boundary condition of  $v$ in \eqref{n1}, we can compute out that
	\begin{equation*}
	\partial_{3}\rho_t+v _{\mm{h}}\cdot \nabla_{\mm{h}} \partial_{3} \rho+\partial_{3}v_3\partial_{3}\rho=0 \text{ on } \partial \Omega.
	\end{equation*}
	Following the argument of \eqref{n22a} with the help of the incompressible condition in \eqref{1a}$_3$ and the initial condition $ \partial_3\rho^0 =0$, we easily derive from  the above  identity that
	\begin{equation}\label{n22sa}
	\partial_3\rho |_{\partial\Omega}=0.
	\end{equation}
	
	Similarly, applying $\partial_3^2$ to \eqref{1a}$_1$, and then making use of the boundary conditions of
	$(\bar{\rho}'',\varrho,v )$  in \eqref{n1}, \eqref{0102n1} and \eqref{n22a}, and the incompressible condition, we have
	\begin{equation*} 	\partial^2_{3}\varrho_t+v_{\mm{h}}\cdot \nabla_{\mm{h}}\partial^2_{3} \varrho +2\partial_{3}v_3\partial^2_{3}\varrho=0 \text{ on } \partial \Omega,
	\end{equation*}
	which obviously also implies that
	\begin{equation}\label{nfdsa22a}
	\partial_3^2\varrho |_{\partial\Omega}=0.
	\end{equation}
	Putting \eqref{n22a}--\eqref{nfdsa22a} together yields \eqref{varrho}.
	
	(2) Applying curl to the momentum equation \eqref{1a}$_2$,  we obtain the following curl equation:
	\begin{align}
	\label{l41a}
	& {\rho}  \omega_t -\kappa\bar{\rho}'(-\partial_2,\partial_1,0)^\top\Delta{\varrho}-
	(g-\kappa\bar{\rho}''')(-\partial_2,\partial_1,0)^\top{\varrho}\nonumber  \\
	&=\mu\Delta \omega -   {\rho} v\cdot  \nabla \omega - {\mathbf{M}} - \mathbf{N} ,
	\end{align}
	where we have defined that
	\begin{align}
	\label{202401120321011}
	\begin{cases}
	{\mathbf{M}}  :=(-\bar{\rho}'\partial_{t}v_2,\bar{\rho}'\partial_{t}v_1,0)^{\top},\ \mathbf{N}:=  \mathbf{N}^{\mm{m}} +\mathbf{N}^{\mm{c}} + \mathbf{N}^{\mm{k}},\\
	\mathbf{N}^{\mm{m}}:=  (\partial_{2} {\varrho}\partial_{t}v_3-\partial_{3} {\varrho}\partial_{t}v_2,\partial_{3} {\varrho}\partial_{t}v_1-\partial_{1} {\varrho}\partial_{t}v_3,\partial_{1} {\varrho}\partial_{t}v_2-\partial_{2} {\varrho}\partial_{t}v_1)^{\top},\\
	\mathbf{N}^{\mm{c}}:=(\partial_{2}( {\rho}v)\cdot  \nabla v_3-\partial_{3}( {\rho}v)\cdot  \nabla v_2,\partial_{3}( {\rho}v)\cdot  \nabla v_1-\partial_{1}( {\rho}v)\cdot  \nabla v_3,\\
	\qquad\quad  \partial_{1}( {\rho}v)\cdot  \nabla v_2-\partial_{2}( {\rho}v)\cdot  \nabla v_1)^{\top},\\ 	\mathbf{N}^{\mm{k}}:=\kappa(\partial_{3}\varrho\partial_{2}\Delta\varrho-\partial_{2}\varrho\partial_{3}\Delta\varrho,\partial_{1}\varrho\partial_{3}\Delta\varrho-\partial_{3}\varrho\partial_{1}\Delta\varrho,\partial_{2}\varrho\partial_{1}\Delta\varrho-\partial_{1}\varrho\partial_{2}\Delta\varrho)^{\top}.
	\end{cases}	\end{align}
	
	It is easy to see from the boundary condition of $v$ in \eqref{n1} that
	\begin{align}
	&\omega_{\mm{h}}=\ 0\text{ on }\partial \Omega \label{omegsfdaa1}
	\end{align}	and
	\begin{align} &\partial_{3}\omega_{3}=\partial_{3}
	(\partial_{1}v_2-\partial_{2}v_1) =0 \text{ on }\partial \Omega.\label{20222309093154}
	\end{align}	
	
	Thanks to the boundary conditions  of $(\varrho,\partial_3\varrho,\partial_2^2\varrho,v )$ in \eqref{n1} and \eqref{varrho}, it follows from the first two equations of the system of the equations in \eqref{l41a} that
	\begin{align}
	\partial^2_{3}\omega_{\mm{h}}=0\text{ on }\partial \Omega.\label{omega1}
	\end{align}
	In addition, by virtue of the boundary conditions of $(v_{\mm{h}}, \partial_3^2\omega_{\mm{h}})$ in  \eqref{n1} and  \eqref{omega1}, and  the incompressible condition, we  further obtain the   boundary condition  \eqref{2022401291254}, which implies that
	\begin{align}\partial^3_{3}\omega_{3}=\partial^3_{3}(\partial_{1}v_2-\partial_{2}v_1)=0 \text{ on } \partial \Omega.\label{202223090s93154}
	\end{align}
	Consequently, putting \eqref{omegsfdaa1}--\eqref{202223090s93154} together yields \eqref{omega}.		This completes the proof.
	\hfill $\Box$
\end{pf}

Now we establish several elliptic estimates and some Poincar\'e's inequalities for $\varrho$, $v$ and $v_t$.
\begin{lem}		\label{201805141072}\begin{enumerate}
		\item[(1)] We have the  elliptic estimates:
		\begin{align}&\| \varrho \|_{1,2}\lesssim \| \Delta \varrho\|_{1,0} ,\label{2fsa024012918221} \\
		&\| \varrho \|_{2(1+j)}\lesssim \| \Delta^{1+j} \varrho\|_0 \mbox{ for }j=0,\ 1,\label{2024012918221} \\
		&\|  v_t \|_2\lesssim \|( v_t , \Delta v_t )\|_0 ,\label{20240129182211x} \\
		&\|  v \|_3\lesssim \|( v , \nabla \Delta v)\|_0 ,\label{20240129182211} \\
		&\|  v \|_{4}\lesssim \| (v ,\nabla \Delta \omega)\|_0.
		\label{202401291822} \end{align}
		\item[(2)] We have  Poincar\'e's inequalities:
		\begin{align}
		&\|    v_i\|_{1}\lesssim \| \nabla v_i\|_{0} \mbox{ for }1\leqslant i\leqslant 3 , \label{im3}
		\\
		& \|    \partial_t v_i\|_{1}\lesssim
		\begin{cases}
		\| \nabla \partial_t v_i\|_{0}+\|\varrho_t\|_0\| v\|_0  &\mbox{for }i=1,\ 2; \\
		\|\partial_{3} \partial_{t}v_3 \|_0 &\mbox{for }i=3.
		\end{cases}\label{pv2}	\end{align}\end{enumerate}
\end{lem}	
\begin{pf}(1) Recalling the boundary condition of $(\omega_{\mm{h}},\partial_{3}\omega_3)$ in  \eqref{omega}, we use both the elliptic estimates in Lemmas \ref{xfsddfs2212} and \ref{xfs05072212} to deduce that
	$$\| \nabla \partial_{\mm{h}}\omega\|_1\lesssim\|\nabla \partial_{\mm{h}} \omega\|_0+ \| \Delta \partial_{\mm{h}} \omega\|_0.$$
	Similarly, we also have
	$$\| \nabla \partial_3\omega\|_1\lesssim\|   \nabla\partial_3 \omega\|_0+ \|  \Delta\partial_3 \omega\|_0.$$
	Putting the above two estimates together yields
	$$\| \nabla^2 \omega\|_1\lesssim\|   \nabla^2\omega\|_0+ \|  \Delta\nabla  \omega\|_0.$$
	It further follows from the above estimate that
	$$\|  \omega\|_3\lesssim\|   v\|_3+ \|  \Delta\nabla  \omega\|_0,$$
	which, together with  the Hodge-type elliptic estimate   \eqref{202005021302} (with the incompressible condition)   and the interpolation inequality \eqref{201807291850},  further implies \eqref{202401291822}.
	
	Similarly, thanks to the boundary conditions of $
	(\varrho,\partial_3^2\varrho,v)$ in \eqref{n1} and \eqref{varrho}, we also exploit both the elliptic estimates in Lemmas \ref{xfsddfs2212} and \ref{xfs05072212} to deduce \eqref{2fsa024012918221}--\eqref{20240129182211}.
	
	(2) Utilizing the mass equation  \eqref{1a}$_1$, it is easy to see from the momentum equation \eqref{0101}$_2$ that
	\begin{align*}
	\partial_t(\rho v)+\operatorname{div}(\rho v\otimes v)+\nabla(P-
	\kappa\rho\Delta \rho-\kappa|\nabla \rho|^2/2) +
	\kappa\mm{div}(\nabla \rho\otimes \nabla \rho)= \mu\Delta{v}-\rho g\mathbf{e}^3.
	\end{align*}
	Integrating the first two equations in the above system (of equations)  on $\Omega$, and then using the integration by parts and the boundary condition of $v$  in \eqref{n1}, we get that
	\begin{align}\label{2024013102110}
	\partial_t\int \rho v_i \mm{d}x =0\mbox{ for }i=1,\ 2,
	\end{align}which, together with the initial null condition	$\int \rho^0 v_i^0\mm{d}x=0$, yields \begin{equation*}
	\int \rho v_i\mm{d}x=0.
	\end{equation*}
	Making use of the above null condition,	the (upper and lower) bounds of density in \eqref{im1a} and \eqref{im1},  and 	the generalized Korn--Poincar\'e inequality \eqref{pro4a12}, we obtain \eqref{im3} for   $i=1$ and $2$.
	
	It follows from \eqref{2024013102110} that
	\begin{align}
	\int\rho \partial_tv_i \mm{d}x = -\int \rho_t v_i \mm{d}x  ,\label{imm2}
	\end{align}
	Exploiting  \eqref{im1a}, \eqref{im1}, \eqref{pro4a12} and H\"older's inequality, we arrive at \eqref{pv2} for $i=1$ and $2$.
	
	In addition, \eqref{im3} and \eqref{pv2} obviously hold  for $i=3$ by the boundary condition of $v_3$ and the Poincar\'e's inequality  \eqref{2202402040948}. This completes the derivation of \eqref{im3} and \eqref{pv2}. \hfill $\Box$
\end{pf}

Now we turn to the derivation of  the stabilizing estimates.
\begin{lem}\label{lem:202012242115}
	\emph{Under the sharp stability condition \eqref{2020102241504} and
		the stabilizing condition \eqref{2022205071434}},  it holds that
	\begin{align}
	&	\|v_3\|^2_{ 1}\lesssim E_{\mm{L}}( v_3) \label{F},  \\
	&\left\|\int_0^t v_3(x,s)\mm{d}s\right\|^2_{ 1}\lesssim E_{\mm{L}}\left( \int_0^t {v}_3(x,s)\mm{d}s\right)  \label{sF},\\&\|\partial^i_{\mm{h}} \varrho\|^2_{1}\lesssim E(\partial^i_{\mm{h}}\varrho) \mbox{ for }1\leqslant i \leqslant 2,  \label{2020401300307}
	\end{align} 	
	see \eqref{edeE} and \eqref{eeE} for the definitions of  $E_{\mm{L}} $ and $E $, resp..
\end{lem}
\begin{pf}
	Since $\int_{(0,2\pi L_1)\times (0,2\pi L_2)}\partial_3(\partial^i_{\mm{h}} \varrho/\bar{\rho}')\mm{d}x_{\mm{h}}=0$ for $1\leqslant i\leqslant 3$, by  the classical existence theory  of Stokes problem mentioned in Remark \ref{2022401100126} for given $x_3\in (0,h)$, there exists a unique function pair  $(\phi,\varphi)\in H^2(2\pi L_1\times 2\pi L_2 )$  such that
	$\int_{2\pi L_1\times 2\pi L_2}\phi\mm{d}x_{\mm{h}}=\int_{2\pi L_1\times 2\pi L_2}\varphi\mm{d}x_{\mm{h}}=0$ and
	$$\Delta_{\mm{h}}(\phi,\varphi)^\top+\nabla_{\mm{h}} \theta=0\mbox{ and }\partial_1\phi+\partial_2\varphi  =-\partial_3 (\partial^i_{\mm{h}} \varrho /\bar{\rho}') $$
	for some $\theta\in H^1(0,2\pi L_1)\times (0,2\pi L_2)$.
	We define that
	$${w} := (\phi,\varphi,\partial^i_{\mm{h}} \varrho/\bar{\rho}'),\ v\mbox{ or }\int_0^t v(x,s)\mm{d}s. $$
	Then we can verify that $ {w}\in H_\sigma$ (by a density argument if necessary).

	By virtue of the definition  of  $\kappa_{\mm{C}}$, it holds that
	\begin{align}
	g\int \bar{\rho}' {w}_{3}^2\mm{d}y\leqslant \kappa_{\mm{C}}   \|\bar{\rho}'\nabla   {w}_{3} \|^2_0\mbox{ for any }w\in H_\sigma  ,\nonumber
	\end{align}
	which,  together with the Poincar\'e's inequality
	\eqref{2202402040948},  the the sharp stability condition \eqref{2020102241504} and the stabilizing condition \eqref{2022205071434},  implies that
	\begin{align}
	\|{w}_{3}\|_1^2\lesssim 		\| \nabla  {w}_{3}\|_{0}^2\lesssim  (\kappa-\kappa_{\mm{C}}){\|\bar{\rho}'\nabla  {w}_{3}\|_0^2}\leqslant\kappa \| \bar{\rho}'\nabla {w}_3\|_0^2 -g\int\bar{\rho}' {w}^2_3\mm{d}y =E_{\mm{L}}(w_3) . \label{2022401311434}
	\end{align}
	Both the  estimates \eqref{F} and \eqref{sF} follow from \eqref{2022401311434} by taking $w_3=v_3$ and $\int_0^t v_3(x,s)\mm{d}s$, resp..
	
	Let $ \Upsilon= \partial^i_{\mm{h}} \varrho /\bar{\rho}'$.
	Exploiting the boundary condition of $\varrho$ and the integration by parts, it is easy to compute out that 	\begin{align}
	&\int(|\nabla \partial^i_{\mm{h}} \varrho|^2+{\bar{\rho}'''}|\partial^i_{\mm{h}} \varrho|^2/\bar{\rho}')\mmd x\nonumber \\&=\int(|\bar{\rho}'\nabla_{\mm{h}}\Upsilon   |^2+(\bar{\rho}''\Upsilon +\bar{\rho}'\partial_{3}\Upsilon )^2 +\bar{\rho}'''\bar{\rho}'\Upsilon ^2)\mmd x\nonumber \\&=\int(|\bar{\rho}'\nabla_{\mm{h}}\Upsilon   |^2+(\bar{\rho}''\Upsilon +\bar{\rho}'\partial_{3}\Upsilon )^2 -(\bar{\rho}'')^2\Upsilon ^2-2\bar{\rho}''\bar{\rho}'\Upsilon \partial_{3}\Upsilon )\mmd x = \int |\bar{\rho}'\nabla  \Upsilon|^2 {\mm{d}x}. \nonumber
	\end{align}	
	We further derive from the above relation and \eqref{2022401311434} with $w_3=\partial^i_{\mm{h}} \varrho/\bar{\rho}' $ that
	\begin{align*}
	&  \|   \partial^i_{\mm{h}} \varrho/\bar{\rho}' \|^2_1  =\|  \Upsilon\|^2_1  \lesssim E_{\mm{L}}( \Upsilon)=  E( \partial^i_{\mm{h}} \varrho) ,
	\end{align*}
	which, together with   the stabilizing condition \eqref{2022205071434}  and Young's inequality, yields \eqref{2020401300307}.
	\hfill $\Box$
\end{pf}

Finally, we derive a basic energy estimate for the CRT problem.
\begin{lem}\label{varrho1}	It holds that
	\begin{align}
	& E_{\mm{L}}\left(  \int_0^tv_3(x,\tau)\mm{d}\tau\right)+\|   v\|^2_{0} +c\int_0^t\|v \|^2_{ 1}\mm{d}\tau\nonumber \\
	&\lesssim \|v^0\|_0^2+\left(\| \varrho^0\|_1+\|v^0\|_0\right)\left\| \int_0^tv_3(x,\tau)\mm{d}\tau\right\|_0
	+ \int_0^t\mathfrak{N}(\tau)\mm{d}\tau,
	\label{lemiq11x}
	\end{align}
	see \eqref{202402021099} for the definition of $\mathfrak{N}(\tau)$.
\end{lem}
\begin{pf}
	Taking the inner product  of \eqref{1a}$_{2}$ and  $  v$    in $L^2$, and then using  the boundary condition of $v_3$ in \eqref{n1}, the incompressible condition in \eqref{1a}$_3$, the
	integration by parts and the mass equations \eqref{1a}$_1$, we  obtain	
	\begin{align}	
	&\frac{1}{2}\frac{\mm{d}}{\mm{d}t}\int {\rho} | v|^2\mm{d}x +\mu\int |\nabla v|^2\mm{d}x
	\nonumber\\&=\int \left(\frac{{\varrho}_t v  }{2}- {\rho} v\cdot  \nabla v -g{\varrho}\mathbf{e}^3 - \kappa( \bar{\rho}''\nabla{\varrho}
	+\bar{\rho}'\Delta{\varrho}\mathbf{e}^3+\Delta{\varrho}\nabla{\varrho})\right)\cdot v \mm{d}x,
	\nonumber\\&=\int \left( g\bar{\rho}'\int_0^t v_3(x,\tau)\mm{d}\tau\mathbf{e}^3+ \kappa\left( \bar{\rho}''\nabla \int_0^t\bar{\rho}'v_3(x,\tau)\mm{d}\tau
	\right.\right.\nonumber \\
	&\quad \left.\left.\left.+\bar{\rho}'\Delta \int_0^t\bar{\rho}'v_3(x,\tau)\mm{d}\tau\mathbf{e}^3\right)
	-\left( g\left(\varrho^0-\int_0^tv\cdot \nabla \varrho \mm{d}\tau\right) \right.\right.\right.
	\nonumber\\
	& \quad\left.
	\left. \left. -\kappa\left( \bar{\rho}'''  \left(\varrho^0-\int_0^tv\cdot \nabla \varrho \mm{d}\tau\right)-\bar{\rho}'\Delta \left(\varrho^0-\int_0^tv\cdot \nabla \varrho \mm{d}\tau\right)\right)\right)\mathbf{e}^3
	-\kappa\Delta{\varrho}\nabla{\varrho}\right)
	\right) \cdot v \mm{d}x\nonumber \\
	&=-\frac{1}{2}\frac{\mm{d}}{\mm{d}t}E_{\mm{L}}\left( \int_0^tv_3(x,\tau)\mm{d}\tau\right)+\mathcal{I}_1(t)
	,
	\label{2022401sdfa211727}
	\end{align}
	see \eqref{202402040453} for the definition of $\mathcal{I}_1(t)$.

	Exploiting H\"older's  inequality, the boundary condition of $v_3$ in \eqref{n1},   the integration by parts  and the product estimate  \eqref{fgestims}, we have
	\begin{align*}
	\int_0^t\mathcal{I}_1(\tau)\mm{d}\tau\lesssim & \left(\| \varrho^0\|_1+\|v^0\|_0\right)\left\| \int_0^tv_3(x,\tau)\mm{d}\tau\right\|_0\\
	&+\int_0^t\left(\left\|\int_0^\tau(v_{\mm{h}}\cdot \nabla_{\mm{h}}\varrho+v_3\partial_3\varrho)(s)\mm{d}s\right\|_1\|v_3\|_1+\|\varrho\|_2\|(v_{\mm{h}}\cdot \nabla_{\mm{h}}\varrho+v_3\partial_3\varrho)\|_0\right)\mm{d}\tau\\
	\lesssim &\left(\| \varrho^0\|_1+\|v^0\|_0\right)\left\| \int_0^tv_3(x,\tau)\mm{d}\tau\right\|_0+\int_0^t\mathfrak{N}(\tau)\mm{d}\tau .
	\end{align*}
	Finally, integrating  \eqref{2022401sdfa211727} over $(0,t)$, and then utilizing  the (upper and lower) bounds bounds of density in \eqref{im1a} and \eqref{im1}, Poincar\'e's inequality \eqref{im3} with $i=3$, we  arrive at \eqref{lemiq11x}.  This completes the proof. \hfill $\Box$
\end{pf}

\subsection{A total energy  inequality}
To begin with, we shall  establish a series of basic energy estimates for the highest-order spacial derivatives and the  temporal derivatives of solutions in the following two lemmas.
\begin{lem}\label{lem1}
	It holds that
	\begin{align}
	& \frac{\mm{d}}{\mm{d}t}\left(\kappa\| \Delta^2\varrho\|^2_{0} +\|\sqrt{{\rho}}\Delta\omega\|^2_{0} -
	8\kappa\int \bar{\rho}^{(4)}
	\partial_3 (\varrho/\bar{\rho}') \partial_3^4\varrho\mm{d}x
	\right)+c\|\nabla \Delta\omega\|^2_{0}\nonumber \\
	& \lesssim (\|\varrho\|_{1,2}+\|v_t\|_1 )\|v\|_{4}+
	\|\varrho\|^2_4\|v\|_{1,2}  +\|v\|_3\|v_t\|_2+ (\sqrt{\mathcal{E}}+\mathcal{E})\mathcal{D},
	\label{lemiq1c}
	\end{align}
	see the definitions of $\mathcal{E} $ and $\mathcal{D}$ in \eqref{di}.
\end{lem}
\begin{pf}
	To begin with, we apply  $\Delta^2$ to the mass equation \eqref{1a}$_1$, resp. the vortex equation \eqref{l41a} to obtain
	\begin{align}
	\Delta^2 (\varrho _t+\bar{\rho}'v_3+v\cdot\nabla \varrho)=0,
	\label{l4asf1a}
	\end{align}
	resp.
	\begin{align}\label{lfsda41a}
	&\Delta ({\rho}\partial_{t} \omega + {\mathbf{M}} -\kappa\bar{\rho}'(-\partial_2,\partial_1,0)^\top\Delta{\varrho}-
	(g-\kappa\bar{\rho}''')(-\partial_2,\partial_1,0)^\top{\varrho})\nonumber  \\
	&=\Delta(\mu\Delta \omega-     {\rho} v\cdot  \nabla \omega -\mathbf{N} ).
	\end{align}
	
	It holds that, by taking the inner products of \eqref{l4asf1a} and $\kappa  \Delta^2\varrho$  in $L^2$,
	\begin{align*}	
	&\frac{1}{2}\frac{\mm{d}}{\mm{d}t}\| \sqrt{\kappa}\Delta^2\varrho  \|^2_{0}
	= -\kappa\int\Delta^2 (\bar{\rho}'v_3+v\cdot\nabla \varrho)\Delta^2\varrho \mm{d}x.
	\end{align*}	
	Taking the inner product  of
	\eqref{lfsda41a}   and  $\Delta  \omega $ in $L^2$, and then using the integration by parts   and the all boundary conditions in \eqref{omega}, we can obtain that	\begin{align*}	
	&\frac{1}{2}\frac{\mm{d}}{\mm{d}t}\| \sqrt{{\rho}} \Delta\omega  \|^2_{0} +\mu\|\nabla \Delta  \omega\|^2_{0}
	\\&=
	\int(\Delta ( \kappa\bar{\rho}' \nabla_{\mm{h}}^{\perp}\Delta{\varrho}+
	(g-\kappa\bar{\rho}''') \nabla_{\mm{h}}^{\perp}{\varrho}-{\mathbf{M}}_{\mm{h}}) \cdot  \Delta\omega_{\mm{h}}  \\&\quad  + \left(\frac{\varrho_t\Delta \omega }{2}-\Delta{\rho}\partial_{t} \omega  -(\nabla{\rho}\cdot \nabla ) \omega_t - \Delta ({\rho} v\cdot  \nabla \omega+\mathbf{N} )\right) \cdot  \Delta\omega )\mm{d}x.
	\end{align*}
	Summing up the above two identities, and then using   the integration by parts, the mass equation and the boundary condition of $v_3$ in \eqref{n1}, we arrive at
	\begin{align}
	&\frac{1}{2}\frac{\mm{d}}{\mm{d}t}\|(\sqrt{\kappa}\Delta^2\varrho ,\sqrt{{\rho}} \Delta\omega )\|^2_{0} +\mu\|\nabla \Delta  \omega\|^2_{0}=\sum_{i=2}^5 \mathcal{I}_i,\label{1e21}
	\end{align}		
	where we have defined that
	\begin{align*}
	\mathcal{I}_2:=& \int(\Delta( (g-\kappa\bar{\rho}''') \nabla_{\mm{h}}^{\perp}{\varrho} - {\mathbf{M}_{\mm{h}}} )\cdot \Delta\omega_{\mm{h}}
	-(\bar{\rho}''\omega_t+ \bar{\rho}'\partial_3\omega_t)\cdot  \Delta\omega  ) \mm{d}x,\\	\mathcal{I}_3:=&
	\kappa\int(\Delta( \bar{\rho}'\nabla_{\mm{h}}^{\perp}\Delta \varrho)\cdot \Delta \omega_{\mm{h}}-\Delta^2  ( \bar{\rho}'v_3) \Delta^2\varrho )\mm{d}x ,\\
	\mathcal{I}_4:=	&  -\int \left( \Delta  \varrho  \omega_t+ (\nabla \varrho  \cdot\nabla)\omega_t+ \partial_i({\rho} v)\cdot  \nabla \partial_i\omega-\Delta(\rho v)\cdot \nabla \omega+ \Delta\mathbf{N}) \right)\cdot  \Delta\omega    \mm{d}x,	 \\
	\mathcal{I}_5:=& -\kappa \int  \Delta^2(v\cdot\nabla {\varrho}) \Delta^2\varrho   \mm{d}x. \end{align*}
	Next we estimate for $\mathcal{I}_2$--$\mathcal{I}_5$ in sequence by five steps.
	
	(1) Exploiting  the boundary conditions of $ \partial_3^{2i} \omega_{\mm{h}}$ in \eqref{omega}, and the integration by parts, it is easy to see that
	\begin{align}
	\mathcal{I}_2=& \int(\partial_i( (g-\kappa\bar{\rho}''') \nabla_{\mm{h}}^{\perp}{\varrho} - {\mathbf{M}_{\mm{h}}} )\cdot \Delta \partial_i\omega_{\mm{h}}-(\bar{\rho}''\omega_t+ \bar{\rho}'\partial_3\omega_t)\cdot  \Delta\omega  )\mm{d}x\nonumber \\
	\lesssim&(\|\varrho\|_{1,1}+\|v_t\|_{1})\|\omega\|_3+\|\omega\|_2\|v_t\|_2. \label{2024102231268}
	\end{align}
	
	(2)
	Utilizing the boundary  conditions  of $(\bar{\rho}'',\omega_{\mm{h}},\partial_3^{2}\omega_{\mm{h}})$ in \eqref{0102n1} and \eqref{omega},  the integration by parts and the relation
	\begin{align}
	\label{20223402201565}
	\Delta v_3=\partial_{2}\omega_1-\partial_{1}\omega_2 \mbox{ (by the incompressible condition)},
	\end{align}  we get that
	\begin{align}	\mathcal{I}_3= &-
	\kappa\int(\partial_i( \bar{\rho}'\nabla^\bot_{\mm{h}}\Delta \varrho)\cdot \partial_i \Delta\omega_{\mm{h}}+\Delta  (\bar{\rho}''' v_3 +2 \bar{\rho}'' \partial_{3} v_3 +\bar{\rho}'\Delta v_3) \Delta^2\varrho)\mm{d}x
	\nonumber  \\
	= &
	\kappa\int( \Delta  ( \bar{\rho}'(\partial_{1}\omega_2-\partial_{2}\omega_1)) \Delta^2\varrho -
	\partial_i ( \bar{\rho}'\nabla^\bot_{\mm{h}}\Delta \varrho)\cdot \partial_i \Delta \omega_{\mm{h}})\mm{d}x + 	\mathcal{I}_{3,1}
	=  \mathcal{I}_{3,1}+ 	\mathcal{I}_{3,2},
	\label{2022401292957}
	\end{align}
	where we have defined that
	\begin{align*}	
	&\mathcal{I}_{3,1 } = - \kappa\int\Delta  (\bar{\rho}''' v_3 +2 \bar{\rho}'' \partial_{3} v_3  )\Delta^2\varrho\mm{d}x,	\\
	&\mathcal{I}_{3,2}:= \kappa\int(( \partial_i \Delta( \bar{\rho}' \omega_{\mm{h}})- \bar{\rho}'  \partial_i  \Delta\omega_{\mm{h}} )\cdot \nabla^\bot_{\mm{h}} \partial_i\Delta\varrho -
	\bar{\rho}''\nabla^\bot_{\mm{h}}\Delta \varrho \cdot \partial_3 \Delta \omega_{\mm{h}})\mm{d}x.
	\end{align*}
	
	Making use of the boundary conditions  of $(\bar{\rho}'',\varrho,\partial_3\varrho,\partial_3^{2}\varrho,v)$ in \eqref{n1}, \eqref{0102n1} and \eqref{varrho},
	H\"older's inequality,  the incompressible condition,  Newton--Leibniz formula and the integration by parts, we  derive that
	\begin{align}	\mathcal{I}_{3,1}
	=&\kappa\int(\nabla\Delta ( \bar{\rho}''' v_3 +2 \bar{\rho}'' \partial_3 v_3)\cdot \nabla \Delta \varrho\mm{d}x  - 4\kappa\int_{\partial\Omega} \bar{\rho}^{(4)}\partial_3 v_3 \partial_3^3\varrho\mathbf{n}_3\mm{d}x  \nonumber \\
	=&\kappa\int(\nabla\Delta ( \bar{\rho}''' v_3 -2 \bar{\rho}'' \mm{div}_{\mm{h}} v_{\mm{h}}   )\cdot \nabla \Delta \varrho -4\partial_3( \bar{\rho}^{(4)}\partial_3 v_3 )\partial_3^3\varrho )\mm{d}x  +\tilde{\mathcal{I}}_{3,1} \nonumber \\
	=& \kappa \int\left(  \Delta^2\left(\bar{\rho}''' \int_0^{x_3} v_{\mm{h}}(x_{\mm{h}},s)\mm{d}s+2 \bar{\rho}''  v_{\mm{h}} \right) \cdot\nabla_{\mm{h}}\Delta\varrho\right.
	\nonumber\\ &\quad \left.- 4\partial_3^2( \bar{\rho}^{(4)} v_{\mm{h}} )\cdot \nabla_{\mm{h}}\partial_3^2\varrho  \right)\mm{d}x
	+\tilde{\mathcal{I}}_{3,1} \leqslant c \|\varrho\|_{1,2} \|v\|_4+\tilde{\mathcal{I}}_{3,1}, \label{2022420200401}
	\end{align}
	where $\mathbf{n}_3$ denotes the third component of the outward unit normal vector $\mathbf{n}$ to
	$\partial\Omega$ and
	\begin{align}
	\label{20224020201102}
	\tilde{\mathcal{I}}_{3,1}:=-4\kappa \int \bar{\rho}^{(4)}\partial_3 v_3 \partial_3^4\varrho\mm{d}x .
	\end{align}
	
	Thanks to the mass equation, it holds that
	$$
	\tilde{\mathcal{I}}_{3,1}:=4 \kappa\frac{\mm{d}}{\mm{d}t}\int \bar{\rho}^{(4)}
	\partial_3 (\varrho/\bar{\rho}') \partial_3^4\varrho\mm{d}x   +4\kappa\overline{\mathcal{I}}_{3,1},
	$$
	where we have defined that
	$$\overline{\mathcal{I}}_{3,1}=\int \bar{\rho}^{(4)}(
	\partial_3 (\varrho /\bar{\rho}') \partial_3^4(\bar{\rho}'v_3+v\cdot \nabla \varrho)  + \partial_3 (v\cdot \nabla \varrho/\bar{\rho}') \partial_3^4\varrho)\mm{d}x.$$
	Similarly to \eqref{2022420200401}, it is obvious that
	\begin{align} \overline{\mathcal{I}}_{3,1}
	=&\int \Bigg( \bar{\rho}^{(4)}\Bigg( \nonumber  \partial_3 (\varrho /\bar{\rho}')
	\left(\sum_{i=1}^4\partial_3^{i} v_{\mm{h}} \cdot \partial_3^{4-i} \nabla_{\mm{h}} \varrho +\sum_{i=1}^3\partial_3^iv_3 \partial_3^{5-i}\varrho\right) \nonumber \\
	& -\partial_3 (\nabla_{\mm{h}} \varrho /\bar{\rho}')  \partial_3^4\left(\bar{\rho}'\int_0^{x_3}v_{\mm{h}}(x_{\mm{h}},s)\mm{d}s\right)- \mm{div}_{\mm{h}}( \partial_3 (\varrho /\bar{\rho}')  v_{\mm{h}})\partial_3^4\varrho \nonumber \\
	& +\partial_3^3v_{\mm{h}} \cdot \nabla_{\mm{h}} \left(\partial_3 \varrho \partial_3 (  \varrho /\bar{\rho}')\right) + \partial_3 \left(  (v_{\mm{h}} \cdot \nabla_{\mm{h}} \varrho+v_3\partial_3\varrho)/\bar{\rho}'\right) \partial_3^4\varrho \Bigg) \nonumber \\
	&- \partial_3(\bar{\rho}^{(4)}  \partial_3 (\varrho /\bar{\rho}')  v_3) \partial_3^{4}\varrho \Bigg)\mm{d}x \nonumber \\
	\lesssim & \|\varrho\|_{1,1}\|v\|_4+\|\varrho\|_4(\|\varrho\|_{1,3} \|v\|_4+\|\varrho\|_4 (\|v\|_{1,0}+\|
	(v_3,\partial_3  v_3)\|_{2})).
	\label{20224020201102sx}\end{align}
	
	In view of \eqref{20224020201102} and \eqref{20224020201102sx}, we deduce from \eqref{2022420200401} that
	\begin{align*}	\mathcal{I}_{3,1}
	\leqslant& 4\kappa\frac{\mm{d}}{\mm{d}t}\int \bar{\rho}^{(4)}
	\partial_3 (\varrho/\bar{\rho}') \partial_3^4\varrho\mm{d}x
	\nonumber \\
	&+c( \|\varrho\|_{1,2} \|v\|_4+\|\varrho\|_4(\|\varrho\|_{1,3} \|v\|_4+\|\varrho\|_4 (\|v\|_{1,0}+\|(v_3,\partial_3  v_3)\|_{2})).   
	\end{align*}
	In addition, exploiting  the boundary condition of $( \varrho,\partial_3^2\varrho)$  in \eqref{varrho} and the integration by parts, it is easy to see that
	\begin{align*}   \mathcal{I}_{3,2}
	= & \kappa\int(\partial_i\Delta(\bar{\rho}'  \Delta \partial_i \omega_{\mm{h}})- \Delta ^2 ( \bar{\rho}' \omega_{\mm{h}}) )\cdot \nabla^\bot_{\mm{h}}  \Delta \varrho -
	\bar{\rho}''\nabla^\bot_{\mm{h}} \Delta \varrho \cdot \partial_3 \Delta \omega_{\mm{h}})\mm{d}x
	\lesssim   \|\varrho\|_{1,2}\|\omega\|_{3}.
	\end{align*}
	Putting the above two estimates into \eqref{2022401292957} yields
	\begin{align}
	\mathcal{I}_3\leqslant &4\kappa\frac{\mm{d}}{\mm{d}t}\int \bar{\rho}^{(4)}
	\partial_3 (\varrho/\bar{\rho}') \partial_3^4\varrho\mm{d}x\nonumber \\
	&+c ( \|\varrho\|_{1,2} \|v\|_4+\|\varrho\|_4(\|\varrho\|_{1,3} \|v\|_4+\|\varrho\|_4 (\|v\|_{1,0}+\|(v_3,\partial_3  v_3)\|_{2})).    
	\label{20224011120226}
	\end{align}
	
	(3) By the boundary condition of $(\varrho,\partial_3\varrho)$ in \eqref{varrho}, we easily check that
	\begin{align}\partial_3(\mathbf{N}_3^{\mm{m}}+\mathbf{N}_3^{\mm{k}})|_{\partial\Omega}=0.
	\label{2022240123120405}
	\end{align}
	By \eqref{2022240123120405}, the boundary conditions of $(\omega_{\mm{h}},\partial_3^2\omega_{\mm{h}})$ in \eqref{omega}, the integration by parts and the product estimates in \eqref{fgestims},  the integral $\mathcal{I}_4 $ can be estimated as follows:
	\begin{align}
	\mathcal{I}_4 \lesssim & c(  (1+\|\varrho\|_2)\|v\|_2\|\omega\|_2+\|\varrho\|_2\|\omega_t\|_1 )\|v\|_4
	+ \int   \nabla (\mathbf{N}^{\mm{m}}+\mathbf{N}^{\mm{k}})   :\nabla  \Delta\omega    \mm{d}x\nonumber \\
	\lesssim & (\|\varrho\|_{1,3}\|\varrho\|_3+\|\varrho\|_{1,2}\|\varrho\|_4+ (1+\|\varrho\|_2)
	\|v\|_3^2+
	\|\varrho\|_2\|v_t\|_2 )\|v\|_4
	.   \label{20224012292110}
	\end{align}
	
	(4) Obviously, it holds that
	\begin{align}
	\mathcal{I}_5= & -\kappa \int  (\Delta^2 v \cdot\nabla {\varrho}
	+  2 \Delta  v \cdot\nabla\Delta {\varrho} +  4(\partial_{i} \Delta v\cdot \nabla \partial_{i}  {\varrho}
	\nonumber \\
	&+
	\partial_{i} \partial_j v\cdot \nabla \partial_{i} \partial_{j}  {\varrho}
	+ \partial_{i}  v\cdot \nabla \partial_{i}\Delta {\varrho}  )\Delta^2 \varrho   \mm{d}x= \kappa  \sum_{1\leqslant i\leqslant 3}\mathcal{I}_{5,i}  ,\label{un31}
	\end{align}
	where we have defined that
	\begin{align*}
	\mathcal{I}_{5,1}:=&  - \int (  (\Delta^2 v \cdot\nabla {\varrho}
	+  2 \Delta  v \cdot\nabla\Delta {\varrho} +  4(\partial_{i} \Delta v\cdot \nabla \partial_{i}  {\varrho}+
	\partial_{i} \partial_j v\cdot \nabla \partial_{i} \partial_{j}  {\varrho}
	\nonumber \\
	&
	+ \partial_{i}  v \cdot \nabla  \partial_{i}\Delta {\varrho}  ) \Delta_{\mm{h}}(\Delta\varrho
	+ \partial_3^2   \varrho)+  ( \Delta^2 v_{\mm{h}} \cdot\nabla_{\mm{h}} {\varrho}
	+  2 (\Delta  v_{\mm{h}} \cdot\nabla_{\mm{h}}\Delta {\varrho}  \\
	&+ \Delta  v_3 \partial_3\Delta_{\mm{h}} {\varrho}) +  4(\partial_{i} \Delta v_{\mm{h}}\cdot \nabla_{\mm{h}} \partial_{i}  {\varrho}+   \nabla_{\mm{h}} \Delta v_3 \cdot  \nabla_{\mm{h}} \partial_3{\varrho}
	\nonumber \\
	& 
	+\partial_{i} \partial_j v_{\mm{h}}\cdot \nabla_{\mm{h}} \partial_{i} \partial_{j}  {\varrho}
	+
	\nabla_{\mm{h}}\partial_i v_3\cdot\nabla_{\mm{h}} \partial_3\partial_i  {\varrho}\\
	& + \partial_{i}  v_{\mm{h}}\cdot \nabla_{\mm{h}} \partial_{i}\Delta {\varrho}
	+ \nabla_{\mm{h}}v_3\nabla_{\mm{h}}\partial_3\Delta {\varrho} ) )  \partial_3^4 \varrho  ) \mm{d}x, \\
	{\mathcal{I}}_{5,2}:= &   -2  \int    (  (\Delta_{\mm{h}} v_3+\partial_3^2 v_3)\partial_3^3 {\varrho} +  2( \Delta\partial_3 v_3\partial_3^2  {\varrho} +
	\nabla \partial_3 v_3 \cdot \nabla \partial_3^2 {\varrho}
	+ \partial_3 v_3\partial_3^2\Delta {\varrho}  )  )\partial_3^4 \varrho  \mm{d}x,
	\nonumber \\  {\mathcal{I}}_{5,3}:= &   -\int   \Delta^2 v_3\partial_3 {\varrho} \partial_3^4 \varrho  \mm{d}x.
	\end{align*}
	It should be noted that we have used the Einstein convention of summation over repeated indices, where $1\leqslant i$, $j\leqslant 3$.
	
	Exploiting the product estimates in \eqref{fgestims}, we get that
	\begin{align*}
	\mathcal{I}_{5,1}+\mathcal{I}_{5,2}\lesssim \|\varrho\|_{1,3}\|\varrho\|_4\|v\|_4+\|\varrho\|_{4}^2(\|\partial_3v_3\|_2+\|v_3\|_{1,2}).
	\end{align*}
	Similarly, by further using the boundary condition of $(v,\partial_3^3v_{\mm{h}})$ in \eqref{n1} and \eqref{2022401291254}, the incompressible condition and the integration by parts, we obtain  that
	\begin{align*}
	{\mathcal{I}}_{5,3}=&  \int   (\partial_3\Delta  \mm{div}v_{\mm{h}}\partial_3 {\varrho} -\Delta_{\mm{h}}\Delta  v_3\partial_3 {\varrho})\partial_3^4 \varrho  \mm{d}x\\
	=&  \int ( ( \nabla _{\mm{h}} \Delta  v_3\cdot \nabla _{\mm{h}}\partial_3 {\varrho}
	-\partial_3\Delta   v_{\mm{h}}\cdot \nabla _{\mm{h}}\partial_3 {\varrho}  )\partial_3^4\varrho\\
	&- \partial_3( \nabla _{\mm{h}} \Delta  v_3\partial_3 {\varrho}
	-\partial_3\Delta   v_{\mm{h}}\partial_3 {\varrho}    )\cdot \nabla _{\mm{h}}\partial_3^3 \varrho )\mm{d}x
	\lesssim\|\varrho\|_{1,3}\|\varrho\|_4   \|v\|_4.
	\end{align*}
	Putting the above two estimates into \eqref{un31} yields
	\begin{align}
	\mathcal{I}_5
	\lesssim  \|\varrho\|_{1,3}\|\varrho\|_4\|v\|_4+\|\varrho\|_4^2(\|\partial_3v_3\|_2+\|v_3\|_{1,2}).
	\label{un3dfas1}
	\end{align}	 
	Finally, inserting \eqref{2024102231268}, \eqref{20224011120226}, \eqref{20224012292110} and \eqref{un3dfas1} into \eqref{1e21}, and then using the incompressible condition  and Poincar\'e's inequality \eqref{2202402040948}, we immediately arrive at \eqref{lemiq1c}. \hfill $\Box$
\end{pf}

\begin{lem}\label{lem2}
	It holds that
	\begin{align}
	&\|\varrho_t\|_i \lesssim
	\begin{cases}
	(1+\|\varrho\|_3)\|v_3\|_i +\|\varrho\|_{1,i}\|v\|_2 &\mbox{for }i=0,\ 1 ;
	\\
	(1+\|\varrho\|_{i+1})\|v_3\|_i +\|\varrho\|_{1,i}\|v\|_i&\mbox{for }i=2,\ 3,
	\end{cases}\label{202221401177321} \\
	&\| v_{t}\|_0 \lesssim  \|\varrho\|_2^2+\|v\|_2 +\|v\|_2^2 \label{20321},\\
	&\|v_t\|_1 \lesssim   (1+\|\varrho\|_3)   ( \| \varrho \|_{3 }+\| \varrho \|_2^2+\|v\|_3 +\|v\|_2^2),\label{20222140117732}
	\\&
	\frac{\mm{d}}{\mm{d}t}(E_{\mm{L}}(v_3)+\|\sqrt{{\rho}} v_t\|^2_{0}) +c \| v_t\|^2_{1} \lesssim  (\sqrt{\mathcal{E} }+\mathcal{E} )\underline{\mathcal{D}}+\|\varrho\|_{1,2}\sqrt{ \underline{\mathcal{D}} \mathcal{D}},\label{1eemt2}		\\
	& \frac{\mm{d}}{\mm{d}t}\left(
	2\kappa \int \Delta ( \bar{\rho}'\Delta{\varrho} )  \Delta v_3\mm{d}x +\mu\|\nabla \Delta v\|_0^2\right)+c\| v_t\|_2^2\nonumber \\
	&\lesssim \|\varrho\|_{1,1}^2+\|  v_3\|_3^2+\|  v_t\|_0^2+ (\sqrt{\mathcal{E} }+\mathcal{E})\mathcal{D} \label{2022sfa2140117732}. \end{align}
\end{lem}
\begin{pf}
	(1) The estimate \eqref{202221401177321} follows  from the mass equation \eqref{1a}$_1$ and the produce estimates in \eqref{fgestims}.
	
	(2) Taking the inner product of the momentum equation \eqref{1a}$_2$ and $v_t$ in $L^2$, and using the boundary condition of $v_3$, the bounds of density in \eqref{im1}, the integration by parts, the incompressible condition and  the product estimate \eqref{fgestims},  we have
	\begin{align*}c\|v_{t}\|_0^2\leqslant  \|\sqrt{\rho}v_{t}\|_0^2=&
	\int (\mu\Delta v -g{\varrho}\mathbf{e}^3-\kappa(\bar{\rho}''\nabla{\varrho}+
	\bar{\rho}'\Delta{\varrho}\mathbf{e}^3+\nabla{\varrho} \Delta{\varrho})-{\rho} v\cdot  \nabla v)\cdot v_t\mm{d}x\\
	\lesssim &(\|(\varrho,v)\|_2 +(\|\varrho\|_3^2+\|v\|_2^2) )\|v_t\|_0,
	\end{align*}
	which, together with Young's inequality, implies
	\begin{align*}
	\|v_{t}\|_0\lesssim \|(\varrho,v)\|_2+\|\varrho\|_3^2+\|v\|_2^2.
	\end{align*}
	
	(3) Applying $\|\cdot\|_0$ to the vortex equation  \eqref{l41a}, and then exploiting   the product estimates in \eqref{fgestims}, we  obtain that
	\begin{align*}
	\| {\rho}\omega_t\|_0 \lesssim (1+\|\varrho\|_3 ) (\| \varrho \|_{1,2} +\|v_t\|_0+\|v \|_2^2)+\|v\|_3.
	\end{align*}
	Thanks to   the Hodge-type elliptic estimate \eqref{202005021302}  and the lower-bound of density in \eqref{im1}, we further derive from the above estimate  that
	\begin{align}
	\label{202240131012024}\| \nabla v_t\|_0 \lesssim   (1+\|\varrho\|_3) (\| \varrho \|_{1,2} +\|v_t\|_0+\|v \|_2^2)+\|v\|_3,
	\end{align}
	which, together with \eqref{20321}, yields \eqref{20222140117732}.
	
	(4) Applying $\partial_{t}$ to  \eqref{1a}$_{2}$, we get
	\begin{align}\label{t1}
	&\partial_t	({\rho}  v_{t}+  {\rho} v\cdot  \nabla v +\nabla  \beta )\nonumber  \\[1mm]
	&=\partial_t(\mu\Delta v-g{\varrho}\mathbf{e}^3-\kappa(\bar{\rho}''\nabla{\varrho}+ \bar{\rho}'
	\Delta{\varrho}\mathbf{e}^3+\nabla{\varrho} \Delta{\varrho})) .
	\end{align}
	Following the argument of \eqref{2022401sdfa211727},   we  derive from \eqref{t1} that
	\begin{align}
	&\frac{1}{2}\frac{\mm{d}}{\mm{d}t}\left(E_{\mm{L}}\left(  v_3(x,\tau) \right)+
	\|\sqrt{\rho}  v_t\|^2_0\right)+\mu\|\nabla v_t\|^2_0
	\nonumber\\&
	=\int \bigg( (g v\cdot \nabla \varrho  +\kappa  \left(\bar{\rho}'\Delta
	(v\cdot \nabla \varrho) - \bar{\rho}'''
	v\cdot \nabla \varrho  )\right)  \mathbf{e}^3
	\nonumber\\
	& \quad -\kappa\partial_t(\Delta{\varrho}\nabla{\varrho})-\partial_t( {\rho} v) \cdot  \nabla v-\varrho_t v_t\bigg) \cdot v_t \mm{d}x.\label{2022401sdsfafa211727}
	\end{align}
	
	By the integration by parts, the product estimates in \eqref{fgestims} and the boundary conditions of $(\varrho,v_3)$, it is obvious that
	\begin{align*}
	\mathcal{I}_6\lesssim &(\|\varrho\|_{1,2}\|v\|_2+ \|\varrho\|_{3}\|v_3\|_2  +\|\varrho\|_{3}\|\varrho_t\|_1+((1+\|\varrho\|_2)\|v\|_1\\
	&+\|\varrho_t\|_2)\|v_t\|_1+\|\varrho_t\|_1(\|v_3\|_2+\|v\|_{1,1})\|v\|_2)   \|v_t\|_1.
	\end{align*}
	Putting  the above estimate into \eqref{2022401sdsfafa211727}, and then using \eqref{pv2}, we get \eqref{1eemt2}.
	
	(5) Applying $\Delta $ to \eqref{1a}$_2$ yields
	\begin{align}
	{\rho}\Delta v_{t}=&\Delta (\mu\Delta v-\nabla  \beta-g{\varrho}\mathbf{e}^3-\kappa(\bar{\rho}''\nabla{\varrho}+
	\bar{\rho}'\Delta{\varrho}\mathbf{e}^3\nonumber \\
	&+\nabla{\varrho} \Delta{\varrho})-   \rho v\cdot  \nabla v)
	- \Delta {\rho} v_{t} -(\nabla {\rho} \cdot \nabla )v_{t}  . 	
	\label{20220020913}
	\end{align}
	In addition, by the incompressible condition and the boundary condition of $(v,\partial^3_{3}v_{\mm{h}})$ in \eqref{n1} and  \eqref{2022401291254}, we have
	\begin{align}
	\label{2022402041403}
	\partial_3 \Delta   v_{\mm{h}} = \Delta   v_3=0 \mbox{ on } 	\partial\Omega.
	\end{align}
	
	Taking the inner product of \eqref{20220020913} and $\Delta v_{t} $
	in $L^2$, and then exploiting the above boundary condition, the incompressible condition, the integration by parts  and the mass equation, we arrive at that
	\begin{align}
	&\frac{\mu  }{2} \frac{\mm{d}}{\mm{d}t}  \|\nabla \Delta v\|_0^2 +\|\sqrt{\rho}\Delta v_t\|_0^2= \sum_{i=7}^9\mathcal{I}_i,
	\label{2024012200334}
	\end{align}
	where we have defined that
	\begin{align*}
	\mathcal{I}_7:=-\kappa\int \Delta(\bar{\rho}'\Delta{\varrho})\Delta \partial_t v_3\mm{d}x,\
	\mathcal{I}_8:=-\int \Delta( g{\varrho} \mathbf{e}^3+  \kappa \bar{\rho}''\nabla{\varrho} +\kappa \nabla{\varrho} \Delta{\varrho} )\cdot \Delta v_t\mm{d}x
	\end{align*}
	and
	\begin{align*}
	&\mathcal{I}_9:=-\int (\Delta(    \rho v\cdot  \nabla v )
	+ \Delta {\varrho} v_{t} +(\nabla \rho \cdot \nabla )v_{t}  )\cdot \Delta v_t\mm{d}x.
	\end{align*}
	
	Thanks to the mass equation and  the boundary condition of $\Delta v_3$ in \eqref{2022402041403}, it holds
	\begin{align*}
	\mathcal{I}_7 =&-\kappa \frac{\mm{d}}{\mm{d}t}  \int \Delta ( \bar{\rho}'\Delta{\varrho} ) \Delta v_3\mm{d}x
	+\kappa\int  \nabla (\bar{\rho}' \Delta(\bar{\rho}v_3+v\cdot\nabla \varrho ))\cdot \nabla \Delta v_3\mm{d}x\nonumber \\
	\leqslant& -\kappa\frac{\mm{d}}{\mm{d}t}  \int \Delta ( \bar{\rho}'\Delta{\varrho} ) \Delta v_3\mm{d}x +c(\|v_3\|_3+\|\varrho\|_4\|v\|_3)\|v_3\|_3.\end{align*}
	
	Exploiting the boundary condition of $\Delta v_3$ in \eqref{2022402041403},  the incompressible condition, the integration by parts and the product estimate \eqref{fgestims}, we have
	\begin{align*}
	\mathcal{I}_8 =&  \int  (\Delta_{\mm{h}}( ( \kappa \bar{\rho}'''-g  ){\varrho} \mathbf{e}^3-\kappa \nabla{\varrho} \Delta{\varrho} )\cdot \Delta v_t- \kappa\partial_3^2(    \nabla_{\mm{h}}{\varrho} \Delta{\varrho} )\cdot \partial_t \Delta v_{\mm{h}}  \\
	&+\partial_3 \nabla_{\mm{h}}(   (\kappa \bar{\rho}''' -g)\varrho -\kappa \partial_3\varrho \Delta{\varrho} ) \cdot (\partial_t \partial_3 \nabla_{\mm{h}} v_3- \partial_t \partial_3  v_{\mm{h}}))\mm{d}x \\
	\lesssim & (\|\varrho\|_{1,1}+ \|\varrho\|_{1,3}\|\varrho\|_4)\|v_t\|_{2}.
	\end{align*}
	
	In addition,
	\begin{align*}
	\mathcal{I}_9\lesssim  ((1+\| \varrho \|_2)\|v\|_2\|v\|_3 +\|\varrho\|_3\|v_{t}\|_1)\|v_t\|_2 .\end{align*}
	Putting the three estimates into \eqref{2024012200334}, and then making use of the lower-bound  of $\rho$ in \eqref{im1}, the elliptic estimate of $v_t$ in  \eqref{20240129182211x} and  Young's inequality, we immediately get \eqref{2022sfa2140117732}.
	\hfill $\square$ \end{pf}
\begin{lem}\label{202240128}
	It holds that
	\begin{align}
	& E(\partial_2\varrho)+E(\partial_1\varrho) \lesssim \|v
	\|_2^2+ \|v_t\|_0^2+{(\sqrt{\mathcal{E}}+\mathcal{E}) {\mathcal{D}}},\label{P1} \\
	&\|\varrho\|_{{1},3}  \lesssim \|\varrho\|_{1,1} + \| v  \|_4 +\| v_t \|_{2} +  (\sqrt{\mathcal{E}}+(\mathcal{E})^{1/4})\sqrt{\mathcal{D}}.\label{r4}
	\end{align} 
\end{lem}
\begin{pf} (1) Taking the inner product of  the vortex equation \eqref{l41a} and $ \nabla_{\mm{h}}^{\perp}\varrho/{\bar{\rho}'}$ in $L^2$, and then using the integration by parts and the boundary condition of $ \varrho $ in \eqref{varrho}, we get that
	\begin{align}
	E(\partial_2\varrho)+E(\partial_1\varrho)
	=\mathcal{I}_{10},\label{LEM2P1}
	\end{align}	
	where we have defined that
	\begin{align}
	\mathcal{I} _{10}:=& \int \frac{1}{\bar{\rho}'}\Bigg(\Bigg(\frac{\mu\bar{\rho}''} {\bar{\rho}'}\partial_{3}\omega_{\mm{h}} -\mathbf{M} _{\mm{h}}
	-{ {\rho}}\partial_{t}\omega_{\mm{h}}-  \mathbf{N}_{\mm{h}} -{\rho}  v\cdot  \nabla \omega_{\mm{h}}\Bigg)
	\cdot\nabla_{\mm{h}}^{\perp}\varrho
	-\mu {\nabla\omega_{\mm{h}} }\cdot\nabla_{\mm{h}}^{\perp}\nabla\varrho
	\Bigg)\mm{d}x .\nonumber
	\end{align}		
	
	It is easily to deduce  that
	$$\mathcal{I} _{10}\lesssim \|\varrho\|_{1,1} (\|v\|_{2}+(1+\|\varrho\|_2)(\|v_t\|_0 + \|v\|_2^2) +\|\varrho\|_{1,1}\|\varrho\|_3+\|\varrho\|_{1,2} \|\varrho\|_2 ).
	$$
	Putting the above estimate into \eqref{LEM2P1}, we obtain \eqref{P1}.
	
	(2) By  the boundary condition of $\varrho$ in \eqref{varrho} and the vortex equation \eqref{l41a},
	we see that $\eta$  satisfies the following Stokes problem:
	\begin{equation}
	\label{Stokesequson}
	\begin{cases}
	\Delta\eta¡¡=  ({\rho}\partial _t \omega_{\mm{h}} -
	(g-\kappa\bar{\rho}''')\eta-\mu\Delta \omega_{\mm{h}} +  {\rho} v\cdot  \nabla \omega_{\mm{h}} +{\mathbf{M}}_{\mm{h}} + \mathbf{N}_{\mm{h}})/\kappa\bar{\rho}'   ,  \\
	\eta  |_{\partial\Omega} =0,
	\end{cases}
	\end{equation}
	where we have defined that $\eta = \nabla^{\perp}\varrho$.
	Applying the elliptic estimate \eqref{xfsdsaf41252} to the above problem \eqref{Stokesequson}, and then using the product estimate \eqref{fgestims}, we get
	\begin{equation}
	\label{xfsddfsf201705141252zz}
	\|\varrho\|_{{1},3}=\|\eta\|_{3}\lesssim
	\|\varrho\|_{ 1,1}+\|v\|_4+(1+\|\varrho\|_2 )(\| v_t \|_{2}+ \|v\|_3^2)+\| \mathbf{N}  \|_{1}.
	\end{equation}
	In addition,
	\begin{align*}
	\|\mathbf{N} \|_{1}\lesssim & \|\varrho\|_3(\|\varrho\|_{1,3}+ \|v_t\|_1)+\|\varrho\|_{1,2}\|\varrho\|_4
	+(1+\|\varrho\|_2)\|v\|_3^2 \lesssim (\sqrt{\mathcal{E}}+\mathcal{E} )\sqrt{\mathcal{D}} .
	\end{align*}
	Inserting the above estimates into \eqref{xfsddfsf201705141252zz} yields \eqref{r4}. \hfill $\Box$
\end{pf}

Now we are in position to building the total energy inequality   for the CRT problem.
\begin{pro}\label{lem3}
	It holds that
	\begin{align}
	\label{omessetsimQ}
	&\sup_{0\leqslant  t\leqslant  T}{ \mathcal{E} (t)   }+\int_0^{T} \mathcal{D} (t) \mm{d}t\nonumber\\&\lesssim \|(\nabla\varrho^0,{v}^0)\|_3^2(1+\|(\nabla\varrho^0,{v}^0)\|_3^2)+\sup_{0\leqslant  t\leqslant  T}{ \mathcal{E}^3 (t)   } + \sup_{0\leqslant  t\leqslant  T}\mathcal{E}(t) \int_0^{T}\left(\|\varrho\|_{1,1} +\|v\|_{1,2} \right)   \mm{d}t\nonumber \\
	&\quad + \sup_{0\leqslant  t\leqslant  T}\sqrt{\mathcal{E}(t)}\left(1+\sup_{0\leqslant  t\leqslant  T}\mathcal{E}^{1/2}(t)\right) \left(\left(\int_0^T ( \|\varrho\|_{1,1} +\|v_3\|_1 ) \mm{d}\tau\right)^2 +\int_0^{T} \mathcal{D} (t)    \mm{d}t\right),
	\end{align}
	where 
	\begin{align}
	\label{omessesdftsimQ}
	\mathcal{E}\lesssim\|(\nabla\varrho,{v})\|_3^2(1+\|(\nabla\varrho,{v})\|_3^4).
	\end{align}
\end{pro}
\begin{pf}
	Exploiting Young's inequality and the elliptic estimate of $\omega $ in \eqref{202401291822}, it follows from \eqref{lemiq1c}, \eqref{1eemt2}	and \eqref{2022sfa2140117732}  that
	\begin{align}
	\frac{\mm{d}}{\mm{d}t}\tilde{\mathcal{E}}(t)+\tilde{\mathcal{D}}(t)  \lesssim  \chi(\|\varrho\|_{{1},2}^2+\|  v\|_3^2)+\|\varrho\|^2_4\|v\|_{1,2}+\chi^2 (\sqrt{\mathcal{E}}+\mathcal{E} )\mathcal{D}  \label{lemisaq1}
	\end{align}
	for sufficiently large positive constant $\chi\geqslant 1$, where we have defined that
	\begin{align*}
	\tilde{\mathcal{E}}(t):= & \kappa\| \Delta^2\varrho\|^2_{0}+\|\sqrt{{\rho}}\Delta\omega\|^2_{0} +\chi\mu\|\nabla \Delta v\|_0^2 +\chi^2(E_{\mm{L}}(v_3) +\|\sqrt{{\rho}} v_t\|^2_{0})  \\
	&+ \kappa \int \left(2\chi\Delta ( \bar{\rho}'\Delta{\varrho} ) \Delta v_3-8 \bar{\rho}^{(4)}\partial_3 (\varrho/\bar{\rho}') \partial_3^4\varrho\right)\mm{d}x
	\end{align*}
	and
	\begin{align*}
	\tilde{\mathcal{D}}(t):=
	\|v\|^2_4+\chi( \chi  \| v_t\|^2_{1}+\|v_t\|_2^2).
	\end{align*}
	Using the dissipative estimates of $\partial_{\mm{h}}\varrho$ in   \eqref{P1} and \eqref{r4},    the interpolation inequality \eqref{201807291850} and  the stabilizing estimate \eqref{2020401300307}, we derive from \eqref{lemisaq1} that, for sufficiently large positive constant $\chi\geqslant 1$,
	\begin{align}
	\frac{\mm{d}}{\mm{d}t}\tilde{\mathcal{E}}(t)+\tilde{\mathcal{D}}(t)  \lesssim  \chi^2 \|  v\|_2^2 +\chi\|v\|^2_3+\|\varrho\|^2_4\|v\|_{1,2}+\chi^2 (\sqrt{\mathcal{E}}+\mathcal{E} )\mathcal{D}  \label{lemiq1}
	\end{align}
	
	In addition, making use of \eqref{202221401177321}--\eqref{20222140117732}, the bounds of  density in \eqref{im1},   the dissipative estimates of $\partial_{\mm{h}}\varrho$ in   \eqref{P1} and \eqref{r4},  the elliptic estimates of $(\varrho,v )$ in \eqref{2024012918221}  and \eqref{20240129182211},  the interpolation inequality \eqref{201807291850}, Poincar\'e's inequality  \eqref{2202402040948}, the stabilizing estimate \eqref{F}  and Young's inequality, we easily obtain that, for sufficiently large positive constant $\chi$,
	\begin{align}
	{\mathcal{E}}(t)\lesssim & \tilde{\mathcal{E}}(t)+ \|\varrho\|_4^2(\|\varrho\|_3^2+\| \varrho \|_2^4+\|v\|_3^2 +\|v\|_2^4)+\| \varrho\|_0^2+\chi^5\|v\|_0^2 \label{202402081929} \\
	\tilde{\mathcal{E}}(t)\lesssim & \chi^2 \|(\nabla \varrho, v)\|_3^2(1+\|(\nabla \varrho, v)\|_3^2)\label{2024020819291}
	\end{align}
	and
	\begin{align*}
	{\mathcal{D}}(t)\lesssim\tilde{\mathcal{D}}(t)+\chi^5\|v\|_1^2
	+{(\sqrt{\mathcal{E}}+\mathcal{E}^2)\mathcal{D}}.
	\end{align*}
	
	Integrating \eqref{lemiq1} over $(0,t)$,  then we can deduce from the resulting inequality,  \eqref{2022204022002023} and \eqref{lemiq11x} by further exploiting  the above three estimates, the interpolation inequality \eqref{201807291850}, Poincar\'e's inequality  \eqref{2202402040948}, the stabilizing estimate \eqref{sF} and Young's inequality that, for some sufficiently large $\chi$,
	\begin{align}
	&{\mathcal{E}}(t)+\int_0^t {\mathcal{D}}(\tau) \mm{d}\tau
	\nonumber \\
	&\lesssim \chi^2 \| (\nabla \varrho^0,v^0)\|_3^2(1+\| (\nabla \varrho^0,v^0)\|_3^2) +\sup_{0\leqslant  t\leqslant  T}{ \mathcal{E}^3 (t)   }+ \chi^5\Bigg( \int_0^t (\sqrt{\mathcal{E}}+\mathcal{E} )\mathcal{D}\mm{d}\tau\nonumber \\
	&\quad + \left(\int_0^t ( \|\varrho\|_{1,1}\|v \|_1+\|\varrho\|_2 \|v_3\|_1) \mm{d}\tau\right)^2 + \int_0^t (\|\varrho\|^2_4\|v\|_{1,2}+\mathfrak{N}(\tau))\mm{d}\tau\Bigg),\label{lemiq1saf}
	\end{align}
	which, together with \eqref{2202402040948} and the incompressible condition, yields \eqref{omessetsimQ}. In addition, the estimate \eqref{omessesdftsimQ} obviously hold by \eqref{202402081929}, \eqref{2024020819291} and \eqref{2202402040948}. This completes the proof. \hfill $\Box$
\end{pf}

\subsection{Tangential energy inequality with decay-in-time}
Before establishing the tangential energy inequality with decay-in-time for the CRT problem, we shall first derive the tangential (derivatives') estimates,
\begin{lem}\label{vardrho1}	It holds that
	\begin{align}
	\frac{\mm{d}}{\mm{d}t}\left(\sum_{i+j=1}E(\partial_1^i\partial_2^j\varrho) +\|\sqrt{\rho} v\|^2_{1,0}  \right) +c\| v\|^2_{1,1}\lesssim(\sqrt{\mathcal{E}}+\mathcal{E}) \underline{\mathcal{D}}+\|\varrho\|_{1,2}\sqrt{\underline{\mathcal{D}}
		\mathcal{D}}\label{lemiqsaf11}
	\end{align}
	and
	\begin{align}
	& \frac{\mm{d}}{\mm{d}t}\|(\sqrt{\kappa}\Delta \varrho,\sqrt{{\rho}}  \omega)\|^2_{1,0}
	+\|\partial_3 v_3\|_2^2+ \|v_3\|^2_{3} + \|v\|^2_{1,2} \nonumber  \\
	&  \lesssim  \|\varrho\|_{2,1}\|v\|_{1,1}+
	\|v_t\|_1^2 + (\sqrt{\mathcal{E}}+\mathcal{E}) \underline{\mathcal{D}}+ \|\varrho\|_{1,2}\sqrt{ \underline{\mathcal{D}}\mathcal{D}},
	\label{2022401290712} \end{align} 	
	see the definitions of  ${\underline{\mathcal{E}}}$ resp. $\underline{\mathcal{D}} $ in \eqref{2022402200520957}  resp. \eqref{2022402011338}.
\end{lem}
\begin{pf} (1)  Applying $\partial_{\mm{h}} $ to the mass equation  \eqref{1a}$_{1}$ and the momentum equation \eqref{1a}$_2$ yields
	\begin{equation}\label{h1}
	\begin{cases}
	\partial_{\mm{h}}   (\varrho _t+\bar{\rho}'v_3 +v\cdot\nabla \varrho)=0,\\	
	\partial_{\mm{h}}  ({\rho} v_{t}+  \rho v\cdot  \nabla v +\nabla  \beta)  \\[1mm]=
	\partial_{\mm{h}} (\mu\Delta v-g{\varrho}\mathbf{e}^3-\kappa(\bar{\rho}''\nabla{\varrho}+
	\bar{\rho}'\Delta{\varrho}\mathbf{e}^3+\nabla{\varrho} \Delta{\varrho})). 		\end{cases}
	\end{equation}
	Taking the inner product of \eqref{h1}$_1$ and $\partial_{\mm{h}}(
	(\kappa{\bar{\rho}'''-g})\varrho/{\bar{\rho}'}-\kappa\Delta\varrho)$ in $L^2$ yields
	\begin{align*}	
	& \frac{1}{2}\frac{\mm{d}}{\mm{d}t}\int\left(  \kappa|\nabla\partial_{\mm{h}}\varrho|^2+(\kappa{\bar{\rho}'''}-g)|\partial_{\mm{h}}\varrho|^2/{\bar{\rho}'}\right)
	\mm{d}x \\
	&= \int \partial_{\mm{h}}( \bar{\rho}'v_3+v\cdot \nabla \varrho)\partial_{\mm{h}}( \kappa\Delta\varrho-(\kappa \bar{\rho}'''-{g})\varrho/{\bar{\rho}'})\mm{d}x.  \end{align*}
	Taking the inner products of  \eqref{h1}$_2$ and  $\partial_{\mm{h}}v$ in $L^2$, and then using  the mass equation, the	integration by parts, the incompressible condition, and the boundary condition of $v_3$ in \eqref{n1}, we can obtain	
	\begin{align*}	
	&\frac{1}{2}\frac{\mm{d}}{\mm{d}t}\int {\rho} | \partial_{\mm{h}}v|^2\mm{d}x +\mu\int |\nabla \partial_{\mm{h}} v|^2\mm{d}x\nonumber\\
	&= \int \left(   \partial_{\mm{h}}{\varrho}  v_{t} -\partial_{\mm{h}}({\rho} v) \cdot  \nabla v-g\partial_{\mm{h}}{\varrho}\mathbf{e}^3 - \kappa\partial_{\mm{h}}( \bar{\rho}''\nabla{\varrho}
	+\bar{\rho}'\Delta{\varrho}\mathbf{e}^3+ \Delta{\varrho}\nabla{\varrho} )\right)\cdot\partial_{\mm{h}} v \mm{d}x.
	\end{align*}
	Adding the above two identities together yields
	\begin{align}	
	&\frac{1}{2}\frac{\mm{d}}{\mm{d}t}(E(\partial_{\mm{h}}\varrho)+\|\sqrt{\rho} \partial_{\mm{h}}v\|^2_{0})+\mu\|\nabla v\|^2_{1,0} =\underline{\mathcal{I}}_1 ,
	\label{1e}
	\end{align}
	where we have defined that
	\begin{align*}
	\underline{\mathcal{I}}_1:=&\int ( (g-\kappa\bar{\rho}''' )\partial_{\mm{h}}(v\cdot\nabla {\varrho})\partial_{\mm{h}}\varrho/{\bar{\rho}'}+\kappa( \Delta{\partial_{\mm{h}}\varrho}  v-\Delta{\varrho}\partial_{\mm{h}}v )\cdot\nabla{\partial_{\mm{h}}\varrho}
	\\
	& + (\partial_{\mm{h}}{\varrho}  v_{t}
	-  \partial_{\mm{h}} ({\rho} v )\cdot  \nabla v)\cdot \partial_{\mm{h}}v  )\mm{d}x.
	\end{align*}
	
	The above integral can be estimated as follows:
	\begin{align*}
	\underline{\mathcal{I}}_{1} \lesssim  &   (\|\varrho\|_{1,1}( \|\varrho\|_3+\|v_t\|_0)  + \|\varrho\|_{1,1}\|v\|_2+(1+\|\varrho\|_2)\|v\|_{1,1} ) \|v\|_1) \|v\|_{1,1}.
	\end{align*}
	Putting the above estimate into \eqref{1e} and then recalling the definition of $\underline{\mathcal{D}}$, we obtain \eqref{lemiq11}.
	
	(2) We apply  $\partial_{\mm{h}}\Delta$ resp. $\partial_{\mm{h}} $ to the vortex equation \eqref{1a}$_1$  resp. \eqref{l41a} to obtain
	\begin{equation}\label{1aaaB}
	\begin{cases}
	\partial_{\mm{h}} \Delta  (\varrho _t+\bar{\rho}'v_3+v\cdot\nabla \varrho)=0, \\
	\partial_{\mm{h}} ({\rho}\partial_{t} \omega +   {\rho} v\cdot  \nabla \omega + {\mathbf{M}} )\\
	=\partial_{\mm{h}} (\mu\Delta \omega +\kappa\bar{\rho}'(-\partial_2,\partial_1,0)^\top\Delta{\varrho}+
	(g-\kappa\bar{\rho}''')(\partial_2,-\partial_1,0)^\top{\varrho}- \mathbf{N} ).
	\end{cases}
	\end{equation}
	Taking the inner product of \eqref{1aaaB}$_1$ resp.  \eqref{1aaaB}$_2$ and $\kappa \partial_{\mm{h}}\Delta\varrho$ resp.  $ \partial_{\mm{h}}\omega $    resp.  in $L^2$,  then making use of the boundary conditions of $(v_3,\partial_3^{2i}\omega_{\mm{h}},\partial_3^{2(1+i)}\omega_3 )$ in \eqref{n1} and \eqref{omega}, the integration by parts and the mass equation, we  obtain that
	\begin{align}	
	&\frac{1}{2}\frac{\mm{d}}{\mm{d}t}(\kappa\|\Delta \partial_{\mm{h}}\varrho\|^2_{0}+\|\sqrt{{\rho}} \partial_{\mm{h}}\omega\|^2_{0} )+\mu\|\nabla  \partial_{\mm{h}}\omega\|^2_{0}
	:=\sum_{i=2}^4\underline{\mathcal{I}}_{ i},\label{1safe21}
	\end{align}		
	where we have defined that
	\begin{align} 
	&\underline{\mathcal{I}}_2:= \int(((\kappa\bar{\rho}'''-g) \partial_{\mm{h}}\nabla_{\mm{h}}^{\perp}{\varrho}
	+ \kappa \bar{\rho}'\partial_{\mm{h}}\nabla^\top \Delta  \varrho-\mathbf{M}_{\mm{h}})  \cdot \partial_{\mm{h}}\omega_{\mm{h}}  -\kappa \partial_{\mm{h}} \Delta(\bar{\rho}' v_3) \partial_{\mm{h}}\Delta\varrho) \mm{d}x,\nonumber\\
	&\underline{\mathcal{I}}_3:= -\int \left(   \partial_{\mm{h}}\varrho \omega_t+  \partial_{\mm{h}}(\varrho v)\cdot \nabla \omega +  \partial_{\mm{h}} \mathbf{N} \right)\cdot \partial_{\mm{h}}\omega\mm{d}x\mbox{ and }\underline{\mathcal{I}}_4:=-\kappa\int \partial_{\mm{h}}\Delta(v\cdot\nabla {\varrho}) \partial_{\mm{h}}\Delta\varrho\mm{d}x. \nonumber
	\end{align}	
	
	Exploiting  \eqref{20223402201565},  the boundary condition of $ \partial_3\varrho$ in \eqref{varrho} and the integration by parts,  $\underline{\mathcal{I}}_2$ can be estimated as follows:
	\begin{align}
	\underline{\mathcal{I}}_2\lesssim & \kappa\int\left(  \bar{\rho}'\partial_{\mm{h}}\nabla^\top \Delta  \varrho  \cdot \partial_{\mm{h}}\omega_{\mm{h}}  - \bar{\rho}'  \partial_{\mm{h}} (\partial_{2}\omega_1-\partial_{1}\omega_2) \partial_{\mm{h}}\Delta\varrho\right) \mm{d}x\nonumber \\
	& -\kappa\int
	\partial_{\mm{h}}(\bar{\rho}''' v_3+2 \bar{\rho}'' \partial_{3} v_3 ) \partial_{\mm{h}}\Delta\varrho  \mm{d}x+ c(\|\varrho\|_{2,0}+\|v_t\|_{1,0} ) \|v\|_{1,1}\nonumber \\
	=& c(\|\varrho\|_{2,1}+\|v_t\|_{1} ) \|v\|_{1,1}   -\kappa\int\nabla
	(\bar{\rho}''' v_3+2 \bar{\rho}''\partial_3 v_3)\cdot \partial_{\mm{h}}\partial_{\mm{h}} \nabla \varrho \mm{d}x \nonumber\\ \lesssim&
	(\|\varrho\|_{2,1}+\|v_t\|_{1} ) \|v\|_{1,1}+\|\varrho\|_{2,1}\|(v_3,\partial_3v_3)\|_1.
	\label{2022401290711}
	\end{align}
	
	Making use of  the boundary condition of $ \partial_3\varrho $ in \eqref{varrho}, the integration by parts and the product estimate \eqref{fgestims}, we have
	\begin{align*}
	\int \partial_3\varrho\partial_i\partial_{\mm{h}}
	\Delta\varrho  \partial_{\mm{h}}\omega_j\mm{d}x=\int\partial_i\partial_{\mm{h}}
	\nabla \varrho \cdot \nabla( \partial_3\varrho \partial_{\mm{h}}\omega_j)\mm{d}x \lesssim \|\varrho\|_{2,1}\|\varrho\|_3\|\omega\|_{1,1},
	\end{align*}
	where $i=1$, $2$ and $1\leqslant j\leqslant 3$.
	Thanks to the above estimate, it is easy to derive that
	\begin{align}
	\underline{\mathcal{I}}_3\leqslant & (\|\varrho\|_{1,1}\|v_t\|_{1}+\|\varrho\|_{2,1}\|\varrho\|_3+\|\varrho\|_{1,2}^2  \nonumber \\
	&+(1+\|\varrho\|_2)\|v\|_{1,1}\|v\|_{2}+
	\|\varrho\|_{1,1}\|v\|_{2}^2)\|\omega\|_{1,1}. \label{1align}	
	\end{align}
	
	It obviously holds that
	\begin{align}
	\underline{\mathcal{I}}_4  = & - \kappa \int
	\Bigg( \Bigg(\sum_{i=1}^3(  \partial_i^2(\partial_{\mm{h}} v \cdot \nabla  {\varrho})
	+  \partial_i^2  v_{\mm{h}}  \cdot\nabla_{\mm{h}} \partial_{\mm{h}} {\varrho}
	+2 \partial_iv_{\mm{h}}  \cdot \nabla_{\mm{h}}\partial_i\partial_{\mm{h}}  {\varrho} + \partial_i^2   v_3 \partial_3\partial_{\mm{h}}{\varrho}\nonumber\\
	&+2\partial_i v_3 \partial_i \partial_{\mm{h}} \partial_3 {\varrho}  )\Bigg)\partial_{\mm{h}}\Delta_{\mm{h}}\varrho  +\sum_{i=1}^2  \partial_i^2(\partial_{\mm{h}} v_{\mm{h}}\cdot \nabla_{\mm{h}}  {\varrho}) +
	\sum_{i=1}^3 ( \partial_i^2 v_{\mm{h}} \cdot\partial_{\mm{h}} \nabla_{\mm{h}}  {\varrho} + 2\partial_i v_{\mm{h}} \cdot \nabla_{\mm{h}}\partial_i \partial_{\mm{h}} {\varrho}
	\nonumber\\
	& +\partial_i^2 v_3  \partial_{\mm{h}} \partial_3 {\varrho}+ 2\partial_i v_3  \partial_i \partial_{\mm{h}} \partial_3 {\varrho}  ) + \partial_3^2 \partial_{\mm{h}} v_{\mm{h}}\cdot \nabla_{\mm{h}}  {\varrho}+ 2\partial_3 \partial_{\mm{h}} v_{\mm{h}}\cdot \nabla_{\mm{h}}  \partial_3{\varrho}+ \partial_{\mm{h}}v _{\mm{h}}\cdot \partial_3^2\nabla_{\mm{h}}  {\varrho}  ) \partial_{\mm{h}} \partial_3^2  \varrho \Bigg)\mm{d}x
	\nonumber \\
	& +\kappa\tilde{\underline{{\mathcal{I}}}}_4\leqslant c(\|\varrho\|_{2,1} ( \|\varrho\|_{1,2}\| v\|_{2}+\|\varrho\|_{3}\| v\|_{1,2})+ \|\varrho\|_{1,2}^2(\|v_3\|_2+\|v\|_{1,2})) +\kappa\tilde{\underline{{\mathcal{I}}}}_4,\label{j12}
	\end{align}	
	where we have defined that
	$$\tilde{\underline{{\mathcal{I}}}}_4:=  \int
	\left(  2\partial_{\mm{h}}\mm{div}_{\mm{h}} v_{\mm{h}}  \partial_3^2{\varrho} + \partial_3  \partial_{\mm{h}} \mm{div}_{\mm{h}} v_{\mm{h}}
	\partial_3{\varrho} -\partial_{\mm{h}}v_3\partial_3^3  {\varrho}-\sum_{i=1}^2  \partial_i^2(
	\partial_{\mm{h}}v_3\partial_3 \varrho ) \right)\partial_{\mm{h}}\partial_3^2\varrho \mm{d}x. $$
	
	Utilizing the boundary conditions of $(\partial_3\varrho, v_3)$ in \eqref{1ab} and \eqref{varrho},  the integration by parts, and the product estimate \eqref{fgestims},   the integral $\tilde{\underline{{\mathcal{I}}}}_4$ can be estimated as follows:
	\begin{align*}
	\tilde{\underline{\mathcal{I}}}_4 :=& \int
	\Bigg  (\partial_3(\partial_{\mm{h}}v_3\partial_3^3{\varrho}- 2\partial_{\mm{h}}\mm{div}_{\mm{h}} v_{\mm{h}}  \partial_3^2{\varrho}  -   \partial_3  \partial_{\mm{h}} \mm{div}_{\mm{h}} v_{\mm{h}}
	\partial_3^2{\varrho} )\partial_{\mm{h}}\partial_3  \varrho
	\\
	&  + \partial_3^2 \mm{div}_{\mm{h}} v_{\mm{h}}
	\partial_{\mm{h}}( \partial_3{\varrho} \partial_{\mm{h}}\partial_3 \varrho )-\sum_{i=1}^2  \partial_3\partial_i (
	\partial_{\mm{h}}v_3\partial_3 \varrho) \partial_i\partial_{\mm{h}}\partial_3 \varrho  ) \Bigg)\mm{d}x\\
	\lesssim & \|\varrho\|_{1,1}  \|\varrho\|_4\| v\|_{1,2}+ \|\varrho\|_{2,1} \|\varrho\|_{3} \|v\|_{1,2} +\|\varrho\|_{1,2}\|\varrho\|_{1,3}\|v\|_{1,2}.
	\end{align*}	
	Putting the above estimate into \eqref{j12} yields
	\begin{align}
	{\underline{\mathcal{I}}}_4 \lesssim &\|\varrho\|_{1,2}(\|\varrho\|_{2,1}\| v\|_3 + \|\varrho\|_{1,2} \|\partial_3v_3\|_2+\|\varrho\|_{1,3}\|v\|_{1,2}) \nonumber \\
	&   +( \|\varrho\|_{1,2}^2+
	\|\varrho\|_{2,1}  \|\varrho\|_{3}+\|\varrho\|_{1,1}\|\varrho\|_4)\| v\|_{1,2}
	. \label{j122}
	\end{align} 
	
	Finally, inserting the three estimates  \eqref{2022401290711}, \eqref{1align} and  \eqref{j122} into \eqref{1safe21}, and then using the Hodge-type elliptic estimate \eqref{202005021302}, the incompressible condition and the Poincar\'e's inequalities of  \eqref{im3} and \eqref{tpar}, we arrive at  \eqref{2022401290712}.  \hfill $\Box$
\end{pf}\begin{lem}\label{202240sfa128}
	It holds that
	\begin{align}
	&\|\varrho \| _{2, 1}^2\lesssim   \|v \|_{1,2}^2 + \|v_t\|_{1}^2+
	{(\sqrt{\mathcal{E}}+\mathcal{E} )\underline{D}} +\|\varrho\|_{1,2}\sqrt{\underline{\mathcal{D}}\mathcal{D}}.
	\label{Pfsda2}
	\end{align}
\end{lem}
\begin{pf}
	Taking the inner product of  \eqref{1aaaB} and $\partial_{\mm{h}}\nabla_{\mm{h}}^{\perp}\varrho/{\bar{\rho}'}$ in $L^2$, and then using  integration by parts and the boundary condition of $\varrho$ in \eqref{varrho}, we get that
	\begin{align}	
	&E( \partial_{\mm{h}}\partial_2\varrho)+E( \partial_{\mm{h}}\partial_1\varrho)
	=\underline{\mathcal{I}} _{5} ,\label{LEM2P1fa}
	\end{align}
	where we have defined that
	\begin{align}
	\underline{\mathcal{I}} _{5}
	:=& \int \frac{1}{\bar{\rho}'}\Bigg( \partial_{\mm{h}}\Bigg(\frac{\mu\bar{\rho}''} {\bar{\rho}'}\partial_{3}\omega_{\mm{h}} -\mathbf{M} _{\mm{h}}
	-{ {\rho}}\partial_{t}\omega_{\mm{h}}-  \mathbf{N}_{\mm{h}}
	\nonumber \\
	& -{\rho}  v\cdot  \nabla \omega_{\mm{h}}\Bigg)
	\cdot \partial_{\mm{h}}\nabla_{\mm{h}}^{\perp}\varrho
	-\mu { \partial_{\mm{h}}\nabla\omega_{\mm{h}} }\cdot \nabla\partial_{\mm{h}}
	\nabla_{\mm{h}}^{\perp}\varrho
	\Bigg)\mm{d}x.\nonumber
	\end{align}
	
	Making use of the integrating by parts and the product estimate \eqref{fgestims}, we easily estimate that
	\begin{align}
	\underline{\mathcal{I}} _5\lesssim&\|\varrho\|_{2,1}(  \|v\|_{1,2} +(1+\|\varrho\|_3)
	\|v_t\|_{1}{+}\|\varrho\|_{2,1}\|\varrho\|_3+\|\varrho\|_{1,2}^2  \nonumber \\
	&+ (1+\|\varrho\|_2)\|v\|_{2}\|v\|_{1,2}+\|\varrho\|_{1,1}\|v\|_{2}^2).
	\end{align}
	Putting the above estimate into \eqref{LEM2P1fa}, and then using Young's inequality, we obtain \eqref{Pfsda2}.
	\hfill$\Box$
\end{pf}

Now we are in the position to building the tangential energy inequality with decay-in-time.
\begin{pro}
	\label{lem2dx}
	It holds that
	\begin{align}
	&\sup_{0\leqslant  t\leqslant  T}(\langle t\rangle^2 {	{\underline{\mathcal{E}}}(t)})+\int_0^{T}{\langle t\rangle^2 		\underline{\mathcal{D}} (t)} \mm{d}t\nonumber
	\\
	&\lesssim \|(\nabla\varrho^0,	{v}^0)\|_3^2(1+\|(\nabla\varrho^0,{v}^0)\|_3^2)\nonumber \\
	&\quad + \int_{0}^{T} (   \langle t\rangle^2((\sqrt{\mathcal{E}}+\mathcal{E}) \underline{\mathcal{D}}+\left(1+\sup_{0\leqslant  t\leqslant  T}(\langle t\rangle^2\|{\varrho}\|_{1,2}^2)\right) {\mathcal{D}} )\mm{d}t.  \label{ed1A}
	\end{align}
\end{pro}
\begin{pf}
	Utilizing Young's inequality, \eqref{202221401177321}  and \eqref{Pfsda2}, we can derive from \eqref{1eemt2}, \eqref{lemiqsaf11} and  \eqref{2022401290712} that
	\begin{align}
	&\frac{\mm{d}}{\mm{d}t}\tilde{\underline{ {\mathcal{E}}}}(t)+c\tilde{\underline{ {\mathcal{D}}}}(t) \lesssim    \chi((\sqrt{\mathcal{E}}+\mathcal{E}) \underline{\mathcal{D}}+
	\|{\varrho}\|_{1,2}\sqrt{  \underline{\mathcal{D}}\mathcal{D}} ),	 \label{lemsafiqsaf11}
	\end{align}
	where $\chi\geqslant 1$ is a sufficiently large constant, and  we have defined that
	$$ \tilde{\underline{{\mathcal{E}}}}(t):=\chi\left(\sum_{i+j=1}E(\partial_1^i\partial_2^j\varrho)  +\|\sqrt{\rho} v\|^2_{1,0}+E_{\mm{L}}(v_3) +\|\sqrt{\rho} v_t \|_0^2\right)+\|(\sqrt{\kappa}\Delta \varrho,\sqrt{{\rho}}  \omega)\|^2_{1,0}
	$$
	and
	$$
	\tilde{\underline{{\mathcal{D}}}}(t):=\|\varrho\|_{2,1}^2 +\chi\| v_t\|^2_{1}  +\|(v_3,\partial_3v_3)\|_2^2+\| v\|^2_{1,2}.
	$$
	Moreover, making use of \eqref{202221401177321} with $i=1$,  \eqref{20321}, the bounds of  density in \eqref{im1a} and \eqref{im1}, the elliptic estimate of $\partial_{\mm{h}}\varrho$ in  \eqref{2fsa024012918221}, the Hodge-type elliptic estimate \eqref{202005021302}, the incompressible condition,  the Poincar\'e's inequalities \eqref{im3} and  \eqref{tpar},  and the stabilizing estimates of \eqref{F} and \eqref{2020401300307}, we have
	\begin{align}
	{\underline{\mathcal{E}}}\lesssim &\tilde{\underline{{\mathcal{E}}}}(t)+
	\|\varrho\|_{3}^2 \|v_3\|_2^2 +\|\varrho\|_{1,2}^2\|v\|_2^2\lesssim  \chi\|(\nabla \varrho, v)\|_3^2(1+\|(\nabla \varrho, v)\|_3^2)
	\label{2022040203041},\\
	{\underline{\mathcal{E}}}\lesssim& \chi(\|\varrho\|_{2,1}^2+\|v\|_{1,1}^2+\|v_t\|_0^2 )+
	\|\varrho\|_{1,2}^2\end{align}
	and
	\begin{align}
	{\underline{\mathcal{D}}}(t)
	\lesssim  \tilde{\underline{{\mathcal{D}}}}(t)+\|\varrho\|_{2}^2\|v_3\|_1^2 +\|\varrho\|_{1,1}^2\|v\|_1^2.\label{20220sdf40203041}
	\end{align}
	
	Now multiplying \eqref{1eemt2} resp. \eqref{lemsafiqsaf11} by $\chi^2\langle t\rangle$  resp. $\langle t\rangle^2$, then adding the two resulting inequalities together  and filially using \eqref{2022040203041}, we arrive at that
	\begin{align}
	&\frac{\mm{d}}{\mm{d}t}(\chi^2 \langle t\rangle(E_{\mm{L}}(v_3)+\|\sqrt{{\rho}} v_t\|^2_{0}) +\langle t\rangle^2\tilde{ {\underline{\mathcal{E}}}}(t))+c  \langle t\rangle^2 \tilde{\underline{\mathcal{D}}} +\chi^2\langle t\rangle\|v_t\|_1^2 \nonumber \\ 
	& \lesssim \chi\langle t\rangle (\|\varrho\|_{1,2}^2+\| v\|^2_{1,1}  ) {+}
	\chi^2( \| v_3  \|^2_0+\|v_t\|_0^2)
	\nonumber \\
	&\quad +   \chi^2 \langle t\rangle^2 ((\sqrt{\mathcal{E}}+\mathcal{E}) \underline{\mathcal{D}}+
	\|{\varrho}\|_{1,2}\sqrt{   \underline{\mathcal{D}}\mathcal{D}} ).  \label{lemiq11}
	\end{align}
	Integrating the above inequality over $(0,T)$ and then exploiting  \eqref{2022040203041}, \eqref{20220sdf40203041}, the interpolation inequality \eqref{201807291850},  the stabilizing estimate \eqref{F} and Young's inequality, we have  \eqref{ed1A}. This completes the proof. \hfill $\Box$
\end{pf}

\section{Proof of Theorem \ref{thm2}}\label{subsec:08a}

To begin with, we shall state a local well-posedness result for the CRT problem.
\begin{pro}\label{202102182115}
	Let $\mu$, $\kappa$ be positive constants, and  $0<\bar{\rho}\in {C^7}[0,h]$. There exists  $T_0>0$ such that, for any $(\varrho^0,v^0)\in  {H}^{4}_{\bar{\rho}} \times {^0_\sigma {{H}}^3_{\mathrm{s}}}$ satisfying a necessary compatibility condition	and the positive lower-bound condition of initial density $0<\inf\limits_{x\in\Omega}\big\{{\rho}^{0}(x)\big\}$, the CRT problem \eqref{1a}--\eqref{n1}  admits a unique local(-in-time) classical solution $(\varrho,v)$ with an associated pressure $\beta$; moreover $(\varrho,v,\beta)\in{\mathfrak{P}} _{T_0}\times { \mathcal{V}_{T_0} }\times C^0([0,T],\underline{H}^2)$ and \begin{equation*}
	0<\inf\limits_{x\in\Omega}\big\{{\rho}^{0}(x)\big\}\leqslant {\rho}(t,x)\leqslant \sup\limits_{x\in\Omega}\big\{{\rho}^{0}(x)\big\}\mbox{ for any }(t,x)\in I_{T_0}\times\Omega,
	\end{equation*}
	where $\rho^0:=\varrho^0+\bar{\rho}$.
\end{pro}
\begin{pf}Since Proposition \ref{202102182115} can be easily proved by the standard iteration method as in \cite{JFJSWWWN,JFJSZYYO}, we omit the trivial proof. \hfill $\Box$
\end{pf}

Due to the \emph{a priori} energy inequalities in Propositions \ref{lem3} and \ref{lem2dx}, we can easily establish the global solvability
in Theorem \ref{thm2}. Next, we briefly describe the proof.

Let $(\varrho^0,v^0)$ satisfy the assumptions in Theorem \ref{thm2}.
By the embedding inequality \eqref{esmmdforinfty},  there exists a constant $\delta_1>0$ such that, if $\|\varrho^0\|_3\leqslant \delta_1$, it holds that
$0<\inf\limits_{x\in\Omega}\big\{{\rho}^{0}(x)\big\}$.  From now on, we choose  $\delta$ in Theorem \ref{thm2} to be less than $\delta_1$.

In view of Proposition \ref{202102182115}, there exists a unique local classical solution $(\varrho,v, \beta)$ to the CRT problem of
\eqref{1a}--\eqref{n1}  with the maximal existence time $T^{\max}$, which satisfies
\begin{itemize}
	\item for any $a\in I_{T^{\max}}$,
	the solution $(\varrho,v, \beta)$ belongs to  $(\varrho,v,\beta)\in {\mathfrak{P}} _{a}\times { \mathcal{V}_{a} }\times C^0([0,a],\underline{H}^2) $;
	\item $\limsup_{t\to T^{\max} }\| v( t)\|_3=\infty$  if  $T^{\max}<\infty$.
\end{itemize}

Moreover, by the regularity of $(\varrho,v,\beta)$, the solution satisfies \eqref{omessetsimQ}, \eqref{omessesdftsimQ}, \eqref{ed1A} and
\begin{equation}
0<\inf\limits_{x\in\Omega}\big\{{\rho}^{0}(x)\big\}\leqslant {\rho}(t,x)\leqslant \sup\limits_{x\in\Omega}\big\{{\rho}^{0}(x)\big\}\mbox{ for any }(t,x)\in I_{T^{\max}}\times\Omega .
\label{202402081718}
\end{equation}
In particular, by Young's inequality, the  Poincar\'e's inequality \eqref{2202402040948} and the incompressible condition,  there exists positive constants $c_1 \geqslant 1$ and $\delta_2\leqslant \delta_1$ such that
\begin{align}
\label{omessestsimQ}
&\sup_{0\leqslant  t\leqslant  T}({ \mathcal{E} (t)   }+(\langle t\rangle^2 {	{\underline{\mathcal{E}}}(t)}))+\int_0^{T} (\mathcal{D} (t) +{\langle t\rangle^2 		\underline{\mathcal{D}} (t)} )\mm{d}t \leqslant {c}_1 \|(\nabla\varrho^0,v^0)\|_3^2  /2 ,
\end{align}
if
$$\sup_{0\leqslant  t\leqslant  T}(\|(\nabla\varrho,v)(t)\|_3^2+\langle t\rangle^2 \|\varrho(t)\|_{1,2}^2)+\int_0^T\mathfrak{D}(t)\mm{d}t\leqslant  \delta_2^2,$$ where the constants $c_1$ and $\delta_2$  depend on  the domain $\Omega$, and the other known physical parameters/functions, and we have defined that
$$\mathfrak{D}(t):=  \langle t\rangle^2(\|\varrho\|_{1,1}^2+\|v\|_{1,2}^2) .$$

Let $\delta\leqslant \delta_2/\sqrt{2c_1}$ and
\begin{align}
T^{*}=&\sup \left\{\tau\in I_{T^{\max}}~ \left|~  \|(\nabla\varrho,v)(t)\|_3^2 +\langle t\rangle^2 \|\varrho(t)\|_{1,2}^2 +\int_0^t\mathfrak{D}(\tau)\mm{d}\tau\leqslant 2 c_1\delta^2 \ \mbox{ for any }
\ t\leqslant\tau \right.\right\}.\nonumber
\end{align}
Then, we easily see that the definition of $T^*$ makes sense by the fact
$$ \|(\nabla\varrho^0,v^0)\|_3^2 +\|\varrho^0\|_{1,2}^2 \leqslant   {c}_1 \|(\nabla\varrho^0,v^0)\|_3^2 \leqslant  c_1\delta^2   .$$ Thus, to show the existence of a global solution, it suffices to verify $T^*=\infty$. We shall prove this by contradiction below.

Assume $T^*<\infty$, then by Proposition \ref{202102182115} and \eqref{202402081718}, we have
\begin{align}T^*\in I_{T^{\max}}. \label{20222022519850}
\end{align}
Noting that
\begin{equation*}
\sup_{0\leqslant  t\leqslant    {T^{*}}} (\|(\nabla\varrho,v)(t) \|_3^2 + \langle t\rangle^2 \|\varrho(t)\|_{1,2}^2) + \int_0^{T^{*}}\mathfrak{D}(\tau)\mm{d}\tau     \leqslant  2 c_1 \delta^2 \leqslant  {\delta^2_2},
\end{equation*}
then, by the assertion in \eqref{omessestsimQ}, we have
\begin{align*}
\sup_{0\leqslant  t\leqslant    {T^{*}}}({ \mathcal{E} (t)   }+(\langle t\rangle^2 {	{\underline{\mathcal{E}}}(t)}))+\int_0^{T^*} (\mathcal{D} (t) +{\langle t\rangle^2 		\underline{\mathcal{D}} (t)} )\mm{d}t  \leqslant c_1
\|(\nabla \varrho^0,v^0)\|_2^2 \leqslant   c_1 \delta^2.
\end{align*}
In particular,
\begin{align}
\sup_{0\leqslant  t\leqslant   {T^{*}}} (\|(\nabla \varrho,v)(t) \|_2^2 +\langle t\rangle^2 \|\varrho(t)\|_{1,2}^2)+  \int_0^{T^*}\mathfrak{D}(t)\mm{d}t    \leqslant c_1 \delta^2 .   \label{2020103261534}
\end{align}

By \eqref{20222022519850}, \eqref{2020103261534} and the strong continuity $(\nabla\varrho,v)\in C^0([0,T^{\max}), H^3)$, we see that
there is a  constant $\tilde{T}\in (T^*,T^{\max})$, such that
\begin{align}
\sup_{0\leqslant  t\leqslant  {\tilde{T} }} \|(\nabla \varrho,v)(t) \|_3^2+\int_0^{\tilde{T} }\mathfrak{D}(\tau)\mm{d}\tau  \leqslant 2 c_1\delta^2  , \nonumber
\end{align}
which contradicts with the definition of $T^*$. Hence, $T^*=\infty$ and thus $T^{\max}=\infty$.
This completes the proof of the existence of a global solution. The uniqueness of the global solution is obvious due to the uniqueness result of local solutions in Proposition \ref{202102182115}.

\appendix
\section{Analysis tools}\label{sec:09}
\renewcommand\thesection{A}
This appendix is devoted to providing some mathematical results, which have been used in previous sections. We should point out that $\Omega$ and the simplified notations appearing in what follows are  as same as these defined in Section \ref{202402081729}. In addition, $a\lesssim b$ still denotes $a\leqslant cb$ where the positive constant $c$ depends on the parameters and the domain in the lemmas in which $c$ appears.

\begin{lem}
	\label{201806171834}
	Embedding inequality ( \cite[Theorems 4.12]{ARAJJFF}): let $D\subset \mathbb{R}^3$ be a domain satisfying the cone condition, and  $2\leqslant p\leqslant 6$, then
	\begin{align} 
	&\label{esmmdforinfty}\|f\|_{C^0(\bar{D})}= \|f\|_{L^\infty(D)}\lesssim\| f\|_{H^2(D)}\mbox{ for any }f\in  H^2(D).
	\end{align}
\end{lem}
\begin{lem}
	Interpolation inequality in $H^j$ (see \cite[5.2 Theorem]{ARAJJFF}): let $D$ be a domain in $\mathbb{R}^n$ satisfying the cone condition, then, for any given $0\leqslant j< i$, 
	\begin{equation}
	\|f\|_{H^j(D)}\lesssim\|f\|_{L^2(D)}^{1-\frac{j}{i}}\|f\|_{H^i(D)}^{\frac{j}{i}}\lesssim \varepsilon^{-j/(i-j)}\|f\|_{L^2(D)} +\varepsilon\|f\|_{H^i(D)} \label{201807291850}
	\end{equation}  for any $f\in H^i(D)$ and  for any $\varepsilon>0$,   
	where the two constants $c $ in \eqref{201807291850} are independent of $\varepsilon$.
\end{lem}
\begin{lem}\label{xfsddfsf20180508} Product estimates of $H^i$  (see Section 4.1 in \cite{JFJSNS}):
	let $D\subset \mathbb{R}^3$ be a   domain  satisfying the cone condition, then
	\begin{align}
	\label{fgestims}
	\|fg\|_{H^i(D)}\lesssim \begin{cases}
	\|f\|_{H^1(\Omega)}\|g\|_{H^1(D)} & \hbox{for }i=0; \\
	\|f\|_{H^i(D)}\|g\|_{H^2(D)} & \hbox{for }0\leqslant i\leqslant 2;  \\
	\|f\|_{H^2(D)}\|g\|_{H^i(D)}+\|f\|_{H^i(D)}\|g\|_{H^2(D)}& \hbox{for }i=3 ,
	\end{cases}
	\end{align}
	if the norms on the right hand of the above inequalities are finite.		
\end{lem} 
\begin{lem}
	A Poincar\'e-type inequality (see \cite[Lemma 10.6]{YnaBVG}):
	it holds that
	\begin{equation}
	\|{f}\|_{0}\lesssim\|{f}\|_{L^2(\partial\Omega)}+\|{\partial_3f}\|_{0}
	\mbox{ for all }f\in H^{1}.
	\label{2202402040948}
	\end{equation}
\end{lem}
\begin{lem}\label{10220830}
	A Poincar\'e's inequality (see \cite[Lemma 1.43]{NASII04}): let $1\leqslant p<\infty$, and $D$ be a bounded Lipchitz domain in $\mathbb{R}^n$ for $n\geqslant 2$ or a finite interval in $\mathbb{R}$. Then
	\begin{equation}
	\label{poincare}
	\|w\|_{L^p(D)}\lesssim \|\nabla w\|_{L^p(D)}+\left|\int_{D}w\mathrm{d}x\right|\mbox{ for any }w\in W^{1,p}(D).
	\end{equation}
\end{lem}
\begin{rem}\label{10220saf830p}
	In particular, by the above Poincar\'e's inequality, it holds that, for any given $i\geqslant  0$,
	\begin{align}
	&\label{tpar}
	\| w\|_{ {1},i}\lesssim \|w\|_{ {2},i}\mbox{ for any }w\in H^{2+i}.
	\end{align}
\end{rem}
\begin{lem}\label{pro4a1}
	A generalized Korn--Poincar\'e inequality (referring to \cite[Lemma A.9]{JFJSZYYO}): let  $D$ be a bounded domain satisfying
	the cone condition in $\mathbb{R}^n$ for $n\geqslant 2$ and $a$, $b$ be constants.  Assume that $p\geqslant 1$,
	\begin{align}
	&\label{pro4a11}
	0 \leqslant \chi,\ 0<a\leqslant\|\chi\|_{L^{1}(D)},\ \|\chi\|_{L^{p}(D)}\leqslant b,
	\end{align}
	then
	\begin{align}
	&\label{pro4a12}
	\left\|u\right\|_{L^{2}(D)}\lesssim\left\|\nabla u\right\|_{L^{2}(D)}+\left|\int_{D}\chi u\mathrm{d}x\right|\mbox{ for any }u\in H^1(D).		\end{align}
\end{lem}
\begin{lem}\label{pro4a}
	A	Hodge-type elliptic estimate  (referring to \cite[Lemma A.4]{ZHAOYUI}):
	let $i\geqslant1$,
	then
	\begin{align}
	&\label{202005021302}
	\|\nabla w\|_{i-1}
	\lesssim\|(\mm{curl}w,\mm{div}w)\|_{i-1}\mbox{ for any  }w\in H^i_{\mm{s}}.
	\end{align}
\end{lem}\begin{lem}  \label{xfsddfs2212}
	An elliptic estimate for the Dirichlet boundary value condition (referring to  \cite[Lemma A.7]{ZHAOYUI}):
	Let $i\geqslant  0$, $ {f}^1\in H^{i} $ and $ {f}^2 \in H^{i+1/2}(\partial \Omega)$ be given, then there exists a unique solution $u\in H^{i+2} $
	solving the problem:
	\begin{equation*}
	\begin{cases}
	\Delta u =  {f}^1&\mbox{in }   \Omega, \\[1mm]
	u= {f}^2  &\mbox{on }\partial   \Omega;
	\end{cases}
	\end{equation*}
	moreover,
	\begin{equation}
	\label{xfsdsaf41252}
	\|u\|_{{i+2} } \lesssim
	\| {f}^1\|_{i }+ | {f}^2 |_{ H^{i+1/2}(\partial \Omega)}.
	\end{equation}
\end{lem} 
\begin{lem}  \label{xfs05072212}
	An elliptic estimate for the Neumann boundary value condition  (referring to \cite[Lemma 4.27]{NASII04}):
	Let $a$ be a positive constant, $i\geqslant  0$ and $ {f} \in H^{i} $,  then there exists a unique solution $u\in H^{i+2} $
	solving the problem:
	\begin{equation*}
	\begin{cases}
	-a \Delta u =  \mm{div}{f}  &\mbox{in }   \Omega, \\[1mm]
	\partial_{{\mathbf{n}}}u= {f}\cdot {\mathbf{n}}  &\mbox{on }\partial   \Omega,
	\end{cases}
	\end{equation*}
	where  {${\mathbf{n}}$} denotes the outward unit normal vector to
	$\partial\Omega$; moreover,
	\begin{equation}
	\label{xfsddfsf20170saf5141252}
	\|\nabla u\|_{ {1+i} } \lesssim
	\| {f}  \|_{ i } + \|\mm{div}{f}  \|_{i }  .
	\end{equation}
\end{lem}
\begin{lem}  \label{xfsddfsf201805072212}
	A	Stokes  estimate (see \cite[Lemma A.8]{WYJTIKCT}):
	Let $i\geqslant  0$, $ {f}^1\in H^{i} $, $ {f}^2\in H^{i+1} $ and $ {f}^3\in H^{i+1/2}(\partial \Omega)$ be given such that
	\begin{equation*}
	\int_{\Omega} {f}^2\mmd x=\int_{\partial\Omega} {f}^3\cdot {\mathbf{n}}\mmd x_{\mm{h}},
	\end{equation*}	 	where  {${\mathbf{n}}$} denotes the outward unit normal vector to
	$\partial\Omega$. There exists a unique solution $u\in H^{i+2} $ and $p\in H^{i+1} $
	solving the Stokes problem:
	\begin{equation*}
	\begin{cases}
	\Delta u+\nabla p=  {f}^1,\ \mm{div}u=  {f}^2 &\mbox{in }   \Omega, \\[1mm]
	u= {f}^3  &\mbox{on }\partial   \Omega;
	\end{cases}
	\end{equation*}
	moreover,
	\begin{equation}
	\label{xfsddfsf201705141252}
	\|u\|_{H^{i+2} }+\|p\|_{H^{i+1} }\lesssim
	\| {f}^1\|_{H^i }+\| {f}^2\|_{H^{i+1} }+| {f}^3|_{H^{i+1/2}(\partial \Omega)} .
	\end{equation}
\end{lem}
\begin{rem}\label{2022401100126}
	We mention that the above result in Lemma \ref{xfsddfsf201805072212} for the  horizontally periodic domain $\Omega$  can be similarly extended to the case of the domain being $2\pi L_1\mathbb{T}\times 2\pi L_2\mathbb{T}$.
\end{rem}
\begin{lem}\label{261asdas567} A Poincar\'e's inequality with optimal constant:  it holds that
	\begin{align}
	&\label{20210asda8261406}\|\varphi_3\|_{0}^2\leqslant  \|\nabla \varphi_3\|_{0}^2/(\pi^2 h^{-2}+L_{\mm{max}}^{-2})\mbox{ for any }\varphi\in H_\sigma,		
	\end{align}
	see \eqref{20224013022226} and \eqref{20223090318252}  for the definitions $H^1_\sigma$ and $L_{\mm{max}}$, resp.; moreover the above constant $(\pi^2h^{-2}+L_{\mm{max}}^{-2})^{-1}$ is optimal. 	
\end{lem}
\begin{pf}
	For a given function $f\in L^2 $, we define the horizontal Fourier expansion coefficient of $f$ via
	\begin{equation}\label{hftx}
	\hat{f}(\xi_{\mm{h}},x_3) = \int_{(0,2\pi L_1)\times (0,2\pi L_2)} f(x_{\mm{h}},x_3) e^{-\mm{i}x_{\mm{h}}\cdot \xi_{\mm{h}}}\mm{d}x_{\mm{h}},
	\end{equation}
	where $\xi_i\in L_i^{-1}\mathbb{Z}$ for $i=1$, $2$.		
	
	Let $\varphi\in H_\sigma$.	Due to $\mm{div}\varphi=0$, we have
	$$ i\xi_1 \widehat{\varphi_1}+i\xi_2 \widehat{\varphi_2}+\partial_3\widehat{\varphi_3}=0.$$
	Taking $(\xi_1,\xi_2)=(0,0)$ in the above identity and then using the   boundary condition $\widehat{\varphi_3}(0,0,0)=0$, we have
	\begin{align}\widehat{\varphi_3}(0,0,x_3)
	=\partial_3\widehat{\varphi_3}(0,0,x_3)=0.  \label{20222101512150}
	\end{align} 	
	In addition, it is well-known that there exists a function $\psi_0\in H^1_0(0,h)$ such that (see Lemma 4.4 and  (4.25) in \cite{JFJSARMA2019})
	\begin{align}
	\frac{\|\psi_0\|_{L^2(0,h)}}{\|\psi'_0\|_{L^2(0,h)}} = \sup_{\psi\in H^1_0(0,h)}\frac{\|\psi\|_{L^2(0,h)}}{\|\psi'\|_{L^2(0,h)}}=\frac{h}{\pi}.
	\label{202220101519822}
	\end{align}
	By  Parseval's theorem  (see \cite[Proposition 3.1.16]{grafakos2008classical}),  \eqref{20222101512150} and \eqref{202220101519822}, we have
	\begin{align}
	\|\nabla \varphi_3\|^2_0 &=   \frac{1}{(4\pi^2L_1L_2)^2}\sum_{\xi_{\mm{h}}\in L_1^{-1}\mathbb{Z}\times L_2^{-1}\mathbb{Z}}
	(|\xi_{\mm{h}}|^2 \|\widehat{\varphi_3}(\xi_{\mm{h}},x_3)\|^2_{L^2(0,h)} + \|\partial_3\widehat{\varphi_3} (\xi_{\mm{h}},x_3)\|^2_{L^2(0,h)}) \nonumber \\
	&\geqslant  \frac{( \pi^2 h^{-2}+L_{\mm{max}}^{-2})}{16\pi^4( {L_1} {L_2})^2}\sum_{\xi_{\mm{h}}\in L_1^{-1}\mathbb{Z}\times L_2^{-1}\mathbb{Z}}
	\|\widehat{\varphi_3}(\xi_{\mm{h}},x_3)\|^2_{L^2(0,h)}= ( \pi^2 h^{-2}+L_{\mm{max}}^{-2}) \| \varphi_3\|_0^2, \nonumber
	\end{align}
	which implies \eqref{20210asda8261406}.
	
	Now we further prove that the constant $(\pi^2h^{-2}+L_{\mm{max}}^{-2})^{-1}$ is optimal. Without loss of generality, it suffices to consider the case $L_{\mm{max}}=L_{1}$.
	Thus we define that
	\begin{equation}\label{20222113192}
	\varphi := (-{L_{1}} \psi_0'(x_3)\cos (  x_1/L_{1}),0,-\psi_0(x_3)\sin (  x_1 /L_{1}) ),
	\end{equation}
	where $\psi_0(0)=\psi_0(h)=0$.
	It is easy to see that $\varphi \in H_\sigma$.
	Since $\psi_0$ satisfies \eqref{202220101519822}, we get
	\begin{align}
	\frac{\|\varphi _3\|_0^2}{\|\nabla \varphi _3\|_{0}^2} =   \frac{ \|\psi_0(x_3) \|_{L^2(0,h)}^2 }{\|\psi_0'(x_3) \|_{L^2(0,h)}^2 + L_1^{-2} \|\psi_0  (x_3) \|_{L^2(0,h)}^2 }
	=   ( \pi^2 h^{-2}+L_1^{-2})^{-1} .  \label{202221131926}
	\end{align}
	This means that the constant $(\pi^2h^{-2}+L_1^{-2})^{-1}$ is optimal, and thus we completes the proof. \hfill$\Box$				\end{pf}
\vspace{4mm}
\noindent\textbf{Acknowledgements.}
The research of Fei Jiang was supported by NSFC (Grant Nos. 12231016 and 12371233) and the Natural Science Foundation of Fujian Province of China (Grant No. 2022J01105), and  the research of Zhipeng Zhang by NSFC (Grant No. 12101305).
\renewcommand\refname{References}
\renewenvironment{thebibliography}[1]{%
	\section*{\refname}
	\list{{\arabic{enumi}}}{\def\makelabel##1{\hss{##1}}\topsep=0mm
		\parsep=0mm
		\partopsep=0mm\itemsep=0mm
		\labelsep=1ex\itemindent=0mm
		\settowidth\labelwidth{\small[#1]}%
		\leftmargin\labelwidth \advance\leftmargin\labelsep
		\advance\leftmargin -\itemindent
		\usecounter{enumi}}\small
	\def\newblock{\ }
	\sloppy\clubpenalty4000\widowpenalty4000
	\sfcode`\.=1000\relax}{\endlist}

\end{CJK*}
\end{document}